\documentclass[12pt]{article}
\usepackage{setspace} 
\usepackage{amsmath}  
\usepackage{amssymb} 
\usepackage{amsthm}  
\usepackage{natbib}  
\usepackage{graphicx}  
\usepackage{multirow} 
\usepackage{bm} 
\usepackage[total={6.35in,8.6in}]{geometry}

\usepackage{color} 
\usepackage{pdflscape} 

\usepackage{sectsty}
\usepackage{amssymb, enumerate}

\def\boldX{\boldsymbol{X}}
\def\boldx{\boldsymbol{x}}
\def\Tstar{T^{\ast}}
\def\boldeta{\boldsymbol{\eta}}
\def\boldtheta{\boldsymbol{\theta}}
\def\boldbeta{\boldsymbol{\beta}}
\def\geta{g_{\boldsymbol{\eta}}}

\def\etaopt{\boldsymbol{\eta}^{\text{opt}}}
\def\ci{\perp\!\!\!\perp}
\def\wetai{w_{\boldsymbol{\eta} i}}
\def\hatwetai{\hat{w}_{\boldeta i}}

\def\hatetaopt{\hat{\boldeta}^{\text{opt}}}
\newcommand{\ccell}[1]{\multicolumn{1}{c}{#1}}

\newtheorem{theorem}{{\bf Theorem}}

\def\nuAi{\boldsymbol{\nu}_{A i}}
\def\nuetai{\boldsymbol{\nu}_{\boldeta i}}

\begin{document}

\title{\large{\textbf{On Estimation of Optimal Treatment Regimes For\\ Maximizing $t$-Year Survival Probability}}}
\author{\bigskip
\large{Runchao Jiang$^1$, Wenbin Lu, Rui Song, and Marie Davidian}\\
\normalsize{North Carolina State University}}
\date{}
\maketitle

\begin{footnotetext}[1]
{\textit{Address for correspondence: Runchao Jiang, Department of Statistics, North Carolina State University, Raleigh, NC 27695, U.S.A. Email: rjiang2@ncsu.edu.}}
\end{footnotetext}

\baselineskip=21pt

\begin{abstract}
A treatment regime is a deterministic function that dictates personalized treatment based on patients' individual prognostic information. There is a fast-growing interest in finding optimal treatment regimes to maximize expected long-term clinical outcomes of patients for complex diseases, such as cancer and AIDS. For many clinical studies with survival time as a primary endpoint, a main goal is to maximize patients's survival probabilities given treatments. In this article, we first propose two nonparametric estimators for survival function of patients following a given treatment regime. Then, we derive the estimation of the optimal treatment regime based on a value-based searching algorithm within a set of treatment regimes indexed by parameters. The asymptotic properties of the proposed estimators for survival probabilities under derived optimal treatment regimes are established under suitable regularity conditions.  Simulations are conducted to evaluate the numerical performance of the proposed estimators under various scenarios. An application to an AIDS clinical trial data is also given to illustrate the methods.
\bigskip

\end{abstract}

\newpage
\section{Introduction}

\begin{figure}[ht]\label{motivating example}
\centering
\includegraphics[width=0.475\textwidth]{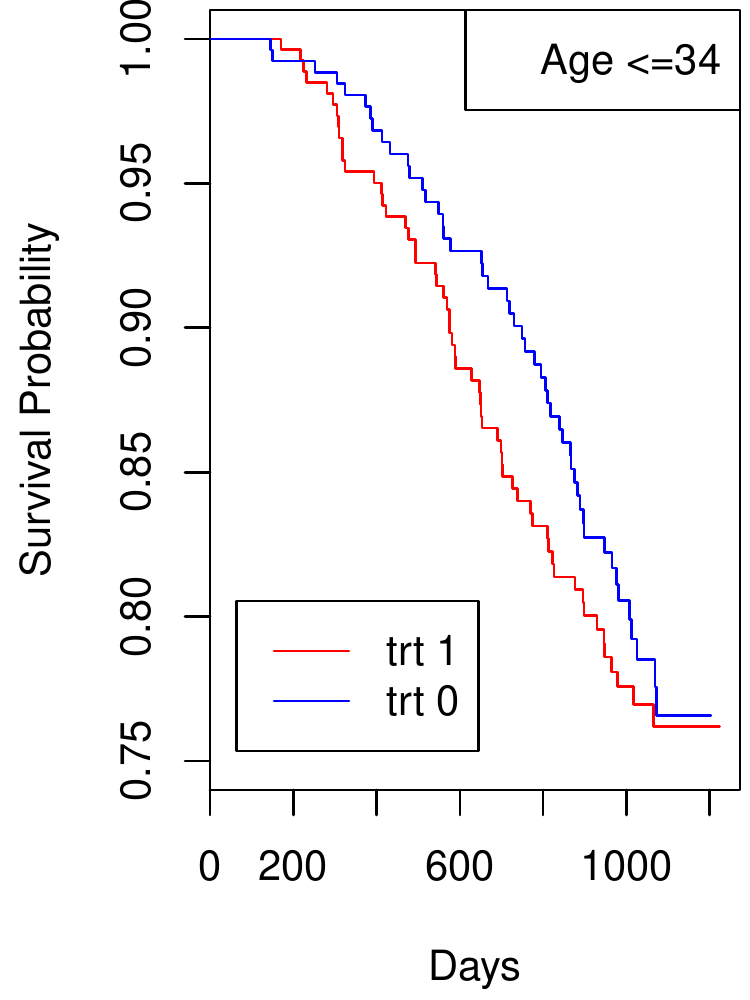}%
\includegraphics[width=0.475\textwidth]{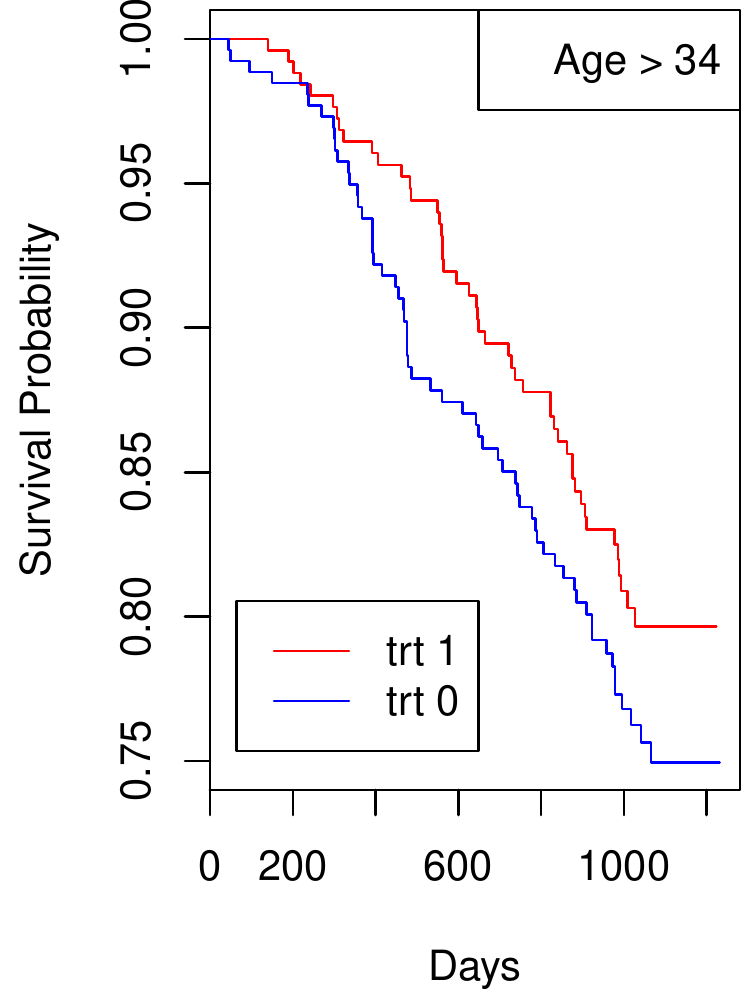}
\caption{Treatment specific Kanplan-Meier curves by age.}
\end{figure}

For many complex diseases, such as cancer, AIDS and mental disorder, there is generally not a uniformly best treatment for all patients. Different patients may favor different treatments, due to individual heterogeneity. For example, in the AIDS Clinical Trials Group Study 175 \citep{Real-Data}, a primary endpoint of interest is the time to having a larger than 50\% decline in the CD4 count, or progressing to AIDS, or death, whichever comes first. We are interested in comparing two treatments: zidovudine plus didanosine (denoted as treatment 1) and zidovudine plus zalcitabine (denoted as treatment 0). We observe that the zidovudine plus zalcitabine treatment is more favorable to younger HIV patients comparing with the zidovudine plus didanosine treatment. To see this, we divide patients into two groups according to the median age of patients, which is 34 in the data. We then plot the treatment specific Kaplan-Meier curves within each age strata, which is given in Figure \ref{motivating example}. From the plot, it can be clearly seen that the zidovudine plus zalcitabine treatment group has almost uniformly larger survival probabilities than the zidovudine plus didanosine treatment group for younger patients with age $\le 34$, while 
the zidovudine plus didanosine treatment group has uniformly larger survival probabilities than the zidovudine plus zalcitabine treatment group for older patients with age $> 34$.

This raises a practically important question on how to appropriately use patients'  individual prognostic information when assigning treatments to maximize an expected long-term clinical outcome of interest, such as $t$-year survival probability. The derivation of optimal individualized treatment regimes, which are a set of treatment decision rules based on patients' individual prognostic information, have received a lot of attention recently, especially for complex diseases such as cancer, AIDS and mental disorder. In addition, for many complex diseases, treatments may be given sequentially at multiple time points. Then a treatment decision rule at a given time point may depend on the baseline prognostic factors, previous assigned treatments and all the intermediate outcomes observed in the past, which results a dynamic treatment regime.  There is a fast development of statistical methods for estimating the optimal dynamic treatment regimes. For example, Q-learning \citep{watkins1989,Q-Learning,murphy2005experimental,Reinforcement} and A-learning \citep{A-Learning,robins2004optimal} are two popular backward induction methods for estimating optimal dynamic treatment regimes. The former is primarily a parametric approach which builds regression models for the so-called Q functions, while the latter is a semiparametric approach which models contrast functions. In addition, A-learning enjoys the double robustness property, i.e. the corresponding estimating equations are asymptotically unbiased when either the baseline mean model or the propensity score model is correctly specified. More recently, \cite{Baqun-Zhang} formularized the problem in a missing data framework and proposed inverse propensity score weighted (IPSW) and augmented IPSW estimators for the expected potential outcome following a specified treatment regime, namely the value function. Then, they proposed to search  the best treatment regime in a pre-specified class of treatment decision rules indexed by parameters to maximize the value function. Such a value-function based optimization method is robust in the sense that it only requires to specify the class of intended treatment regimes but not the models for the Q-functions or contrast functions.  In addition, \cite{OWL-Learning} recast the estimation method of \cite{Baqun-Zhang} in a classification framework and proposed an outcome-weighted learning method to estimate the optimal treatment regime by outcome weighted support vector machines. \cite{Baqun-Zhang2013} extended the value-function based optimization method to estimate the optimal dynamic treatment regime, mainly for two treatment decision time points.

When the outcome of interest is survival time as seen in many clinical trails or observational studies, there is less development for estimation of optimal treatment regimes to maximize patients' survival probabilities given treatments. To our best knowledge, most literatures are focusing on comparing two given treatment regimes. Based on observational experiments with imbalanced treatment assignment, \citet{Chen-and-Tsiatis} and \citet{Min-Zhang} compare the restricted mean survival time for two simple regimes, either giving everyone treatment $1$ or giving everyone treatment $0$.  In addition, \citet{Xiaofei-Bai} proposed doubly-robust estimators for treatment-specific survival probabilities based on observational data with stratified sampling. On the other hand, \citet{Tianxi-Cai} make use of patients' baseline information to predict their risk levels of developing the event of interest at a pre-specified time, i.e. $t$-year survival. Then based on the predicted risk levels, patients are recommended for different therapies accordingly. 
However, this generally can not lead to an optimal treatment regime that maximizes patients' $t$-year survival probabilities. Most recently, \cite{GoldbergKosorok2012} developed a Q-learning algorithm for censored survival data for estimating optimal dynamic treatment regimes and derived its associated finite sample bounds on the generalization error of the policy learned by the algorithm. This approach requires to build a proper regression model for survival times that incorporates both the baseline covariate effects and treatment-covariate interaction effects, which may not be easy in practical applications.

In this article, we propose a value-function based policy search method to estimate the optimal treatment regime that leads to the maximal $t$-year survival probability. Specifically, we first develop two Kaplan-Meier-type estimators for the survival function of patients following a given treatment regime. Then we search the best treatment regime within a class of specified regimes to maximize the associated $t$-year survival probability. Since the estimated $t$-year survival probability following a given treatment regime is a very discrete function of parameters, the direct maximization may be challenging and the resulting estimators may suffer from the numerical instability. To improve the finite sample performance of the estimators, we introduce the kernel smoothing technique to smooth the value function at a proper rate. Both numerical and theoretical properties of the proposed estimators for the $t$-year survival probability following the estimated optimal treatment regime are investigated. In addition, we generalize the proposed method to estimating optimal dynamic treatment regimes and use the case with two treatment decision time points as an illustration.

The rest of the article is organized as follows. We describe our methodology for estimating optimal treatment regimes with a single decision point and multiple decision points in Section~\ref{sec:methodology 1} and \ref{sec:methodology 2}, respectively. The asymptotic properties of the proposed estimators are given in Section~\ref{sec:asymptotic}. Section~\ref{sec:simulation} studies the finite sample performance of the proposed estimators. Section~\ref{sec:application} considers an application to a dataset from the AIDS Clinical Trials Group Study 175 to further illustrate our method. We conclude our work with some discussions in Section~\ref{sec:discussion}. All the proofs are delegated to the Appendix.

\section{Estimation of Optimal Treatment Regime for a Single Decision Time Point} \label{sec:methodology 1}
\subsection{Notation and Assumption}
Consider a study with two treatment options $\mathcal{A}=\{0,1\}$ given at the baseline. For the $i$th patient, $i = 1,\cdots, n$, let $\boldX_i$ denote the $p\text{-dimensional}$ vector of baseline covariates and $A_i$ denote the actual treatment received by the patient. 
In addition, let $T_i$ be the associated continuous survival time of interest, with conditional survival function $S_T(t | a, \boldx) \equiv P (T_i > t | A_i=a, \boldX_i=\boldx)$ and the corresponding conditional cumulative hazard function denoted  by $\Lambda_T(t | a, \boldx)$, where $a = 0/1$. Let $C_i$ denote the right censoring time for patient $i$. The observed data for $n$ independently and identically distributed patients consist of $\{(\boldX_i, A_i, \tilde{T}_i, \delta_i), i=1, \dots, n\}$, where $\tilde{T}_i=\min\{T_i, C_i\}$ and $\delta_i=I\{T_i \leq C_i\}$. Furthermore, we also observe the counting process $N_i(t)=I(\tilde{T}_i \leq t, \delta_i=1)$ and the at risk process $Y_i(t)=I(\tilde{T}_i \geq t)$.

A treatment regime is a deterministic function that maps $\boldX$ to $\mathcal{A}$. For simplicity, we assume the regimes of interest are from $\mathcal{G} = \{\geta: \geta(\boldX)=I\{\boldeta^T \tilde{\boldX} \geq 0\}, \boldeta \in \mathbb{R}^{p+1} , ||\boldeta||=1\}$, where $\tilde{\boldX}=(1, \boldX^{T})^{T}$. However, the proposed method also applies to any other $\mathcal{G}$ that can be indexed by finite-dimensional parameters. Denote the potential survival time of a patient if he/she were given treatment $a$, which may be contrary to fact, as $\Tstar(a)$. Accordingly, define the potential counting process $N^{\ast}(a; t)$ and at risk process $Y^{\ast}(a; t)$ under treatment $a$, where $N^{\ast}(a; t) = I\{ \min (\Tstar(a), C) \leq t, \Tstar(a) \leq C\}$ and $Y^{\ast}(a; t) = I\{ \min (\Tstar(a), C) \geq t\}$.  If a patient follows a given regime $\geta$, we can write the corresponding potential survival time as $\Tstar(\geta) = \Tstar(1)  \geta +\Tstar(0) (1-\geta) $, whose survival function is given by $S^{\ast}(t; \boldeta)=E_{\boldX}[P\{\Tstar(\geta(\boldX)) > t | \boldX\}]$, as well as the potential counting process $N^{\ast}(\geta; t) = N^{\ast}(1; t) \geta + N^{\ast}(0; t) (1-\geta)$ and the potential at risk process $Y^{\ast}(\geta; t) = Y^{\ast}(1; t) \geta + Y^{\ast}(0; t)(1-\geta)$. We are interested in finding the optimal treatment regime in $\mathcal{G}$ that maximizes $t$-year survival probability, that is $\geta^{\text{opt}}(\boldx) \equiv g(\boldx; \etaopt)$, where $\etaopt = \text{arg} \max_{||\boldeta||=1} S^{\ast}(t; \boldeta)$.  Here $t$ is a pre-determined time point, such as 3-year.

To find the optimal treatment regime, we first derive consistent estimators of $S^{\ast}(u; \boldeta)$ for any $u$. To do this, we make the following uninformative censoring assumption: $C$ is independent of $\{T^*(1),T^*(0)\}$ given $A$ and $\boldX$. Let $S_C(t|a, , \boldx) $ denote the survival function of the censoring time given $A=a$ and $\boldX =  \boldx$. If we were able to observe the $\geta \text{-specified}$ potential counting processes $N_i^{\ast}(\geta; s)$'s and at risk processes $Y_i^{\ast}(\geta; s)$'s, an intuitive estimator for $S^{\ast}(u; \boldeta)$ is to consider an inverse probability censoring weighted Kaplan-Meier estimator, specifically,  
\begin{equation}\label{eq:KM}
\widehat{S}^{\ast}(u; \boldeta) = \prod_{s \leq u} \left ( 1-
\frac{ \sum_{i=1}^{n} [d N_i^{\ast}\{\geta(\boldX_i); s\}/S_C\{s|\geta(\boldX_i), \boldX_i\}]}{\sum_{i=1}^{n} [Y_i^{\ast}\{\geta(\boldX_i); s\}/S_C\{s|\geta(\boldX_i), \boldX_i\}]}
\right ).
\end{equation}
However, since $N_i^{\ast}(\geta; s)$'s and $Y_i^{\ast}(\geta; s)$'s are generally  not observable, $\widehat{S}^{\ast}(u; \boldeta)$ is not computable based on observed data. To obtain proper estimators that are computable based on observed data, we make the following two assumptions that are widely used in the causal inference literature \citep{Assumptions}: (i) stable unit treatment value assumption (SUTVA), i.e. $T = T^*(1)A + T^*(0)(1-A)$, and (ii) no unmeasured confounders assumptions, i.e. $\{\Tstar(1), \Tstar(0)\} \ci A | \boldX$. 

\subsection{Estimation Procedure}
Following Zhang et al. (2012), we cast the estimation of $S^{\ast}(u; \boldeta)$ in a missing data framework. Specifically, due to SUTVA, for those patients whose actually received treatment matches with the assigned treatment given by the regime $\geta$,  $N_i^{\ast}(\geta; s) = N_i(s)$ and $Y_i^{\ast}(\geta; s) = Y_i(s)$, which are observed.  For other patients, they are missing. This motivates us to modify the estimator given in (\ref{eq:KM}) by incorporating inverse propensity score weighting. Formally, the weight for the $i\text{th}$ patient is given by
\begin{equation}\label{eq:w}
\wetai = \frac{I[A_i = I\{\boldeta^T\tilde{ \boldX} \geq 0\}] }{\pi(\boldX_i) A_i + \{1-\pi(\boldX_i)\}(1-A_i)}=\frac{A_i I(\boldeta^T\tilde{ \boldX} \geq 0)+(1-A_i)\{1-I(\boldeta^T\tilde{ \boldX} \geq 0)\} }{\pi(\boldX_i) A_i + \{1-\pi(\boldX_i)\}(1-A_i)},
\end{equation}
where $\pi(\boldX_i)=P(A_i=1 | \boldX_i)$ is the propensity score. In practice, $\pi(\boldX_i)$ is either known by design as in randomized clinical trials or needs to be estimated from the data as in observational studies. For the latter case, a parametric model, say a logistic regression is usually used for estimating $\pi(\boldX_i)$, specifically, %
\begin{equation}\label{eq:logistic}
\text{logit}\{{\pi(\boldX_i;\boldtheta)}\} = \boldtheta^{T} \tilde{\boldX_i},
\end{equation}
where $\text{logit}(z)=\log\{z/(1-z)\}$. Let $\hat{\boldtheta}$ denote the maximum likelihood estimator of $\boldtheta$ and define $\hat{\pi}(\boldX_i) = \exp(\hat{\boldtheta}^T \tilde{\boldX}_i)/\{1+\exp(\hat{\boldtheta}^T \tilde{\boldX}_i)\}$. It is known that if the logistic regression model is correctly specified, $\hat{\boldtheta}$ is a consistent estimator of $\boldtheta$. 

To derive the estimator for $S^{\ast}(u; \boldeta)$, we also need to estimate the censoring time survival function $S_C(s|A_i, \boldX_i)$. In many clinical studies with well follow-up, it is reasonable to assume that censoring times are independent of treatment assignment and covariates, i.e. independent censoring assumption. Then, we can use Kaplan-Meier estimator for censoring times to consistently estimate $S_C(s|A_i, \boldX_i)$. For some applications, independent censoring assumption may be restrictive. It can be relaxed to a certain extent. For example, if censoring times are assumed to only depend on treatment assignment, we can use stratified Kaplan-Meier estimators to estimate the treatment-specific censoring time survival function. For more general dependence, we can build a semiparametric model, say a proportional hazards model for censoring times and obtain the model based estimator of $S_C(s|A_i, \boldX_i)$.  For simplicity, from now on we make the independent censoring assumption and let $\hat{S}_C(\cdot)$ denote the Kaplan-Meier estimator for censoring times.    

Let ${\hatwetai}$ denote the estimator of $\wetai$, which is obtained by replacing $\pi(\boldX_i)$ with $\hat{\pi}(\boldX_i)$ in $\wetai$. We propose the following inverse propensity score weighted Kaplan-Meier estimator (IPSWKME) for $S^{\ast}(u; \boldeta)$:
\begin{equation}\label{eq:IPWEK}
\widehat{S}_I(u; \boldeta) = \prod_{s \leq u}
\left\{ 1- \frac{ \sum_{i=1}^n \hatwetai d N_i(s)}{\sum_{i=1}^n \hatwetai Y_i(s)} \right \}.
\end{equation}
Note that the IPSWKME actually dose not depend on the Kaplan-Meier estimator $\hat{S}_C(\cdot)$ for censoring times since it is cancelled out from numerator and denominator under the independent censoring assumption. In Section~\ref{sec:asymptotic}, we will show that $\widehat{S}_I(u; \boldeta)$ is a consistent estimator of $S^{\ast}(u; \boldeta)$ under certain conditions.  Based on $\widehat{S}_I(u; \boldeta)$, the estimated optimal treatment regime to maximize $t$-year survival probability is given by $g(\boldx; \hat{\boldeta}_I^ {\text{opt}})$, where $\hat{\boldeta}_I^ {\text{opt}} = \arg\max_{||\boldeta||=1}\widehat{S}_I(t; \boldeta)$.

Note that the IPSWKME relies on the correct specification of the propensity score model. If it is misspecified, the IPSWKME is generally biased. To improve the robustness of the IPSWKME, we next propose augmented IPSWKME (AIPSWKME) by incorporating assumed model information.
For example, we may posit a proportional hazards (PH) model \citep{PH1} for the conditional cumulative hazard function of $T$ by
\begin{equation}\label{eq:PH}
\Lambda_T (t | A, \boldX) = \Lambda_0(t)\exp \{\boldbeta^T(\boldX^T, A, A\boldX^T)^T \},
\end{equation}
where $\Lambda_0(t)$ is the baseline cumulative hazard function and $\boldbeta$ is a $(2p+1)\text{-dimentional}$ parameter. 
The augmented term for $\wetai d N_i^{\ast}\{\geta(\boldX_i); s\}$ is 
\begin{eqnarray*}
 & &\wetai d N_i^{\ast}\{\geta(\boldX_i); s\} + (1-\wetai) E[d N_i^{\ast}\{\geta(\boldX_i); s\} |\boldX_i]\\
&=&\wetai d N_i^{\ast}\{\geta(\boldX_i); s\} + (1-\wetai) S_T( s | \geta(\boldX_i), \boldX_i) S_C(s)d\Lambda_T(s|\geta(\boldX_i), \boldX_i),
\end{eqnarray*}
where $S_T(s|A_i,\boldX_i)$ and $S_C(s)$ are the conditional survival functions of $T$ and $C$, respectively. Similarly, the augmented term for $\wetai  Y_i^{\ast}\{\geta(\boldX_i); s\}$ is given by $\wetai  Y_i^{\ast}\{\geta(\boldX_i); s\} + (1-\wetai)S_T( s | \geta(\boldX_i), \boldX_i) S_C(s)$. It can be easily shown that the above two augmented terms have the so-called doubly robust property, i.e. they are unbiased when either the propensity score model or the posited PH model is correctly specified. Therefore, we propose the AIPSWKME for $S^{\ast}(u; \boldeta)$ as 
\begin{align}\label{eq:AIPWEK}
\widehat{S}_{A} &(u; \boldeta) \nonumber \\
& = \prod_{s \leq u}
\left ( 1-
\frac{\sum_{i=1}^n [\hatwetai d N_i(s) + (1-\hatwetai) \hat{S}_T\{s | \geta(\boldX_i), \boldX_i\} \hat{S}_C(s) d\hat{\Lambda}_T \{s | \geta(\boldX_i), \boldX_i\}]}
{\sum_{i=1}^n [\hatwetai Y_i(s) + (1-\hatwetai) \hat{S}_T\{s | \geta(\boldX_i), \boldX_i\} \hat{S}_C(s)]}
\right ),
\end{align}
where $\hat{S}_T(s|A_i,\boldX_i)$ is the estimated survival function of $T$ based on the fitted PH model and $\hat{S}_C(s)$ is the Kaplan-Meier estimator for censoring times. 
Based on $\widehat{S}_A(u; \boldeta)$, the estimated optimal treatment regime to maximize $t$-year survival probability is given by $g(\boldx; \hat{\boldeta}_A^ {\text{opt}})$, where $\hat{\boldeta}_A^ {\text{opt}} = \arg\max_{||\boldeta||=1}\widehat{S}_A(t; \boldeta)$. The asymptotic properties of $\widehat{S}_{A} (u; \boldeta)$ and $\widehat{S}_{A} (t; \hat{\boldeta}_A^ {\text{opt}})$ will be studied in Section~\ref{sec:asymptotic}.

\subsection{Computational Aspects}
Note that $\widehat{S}_{I} (t; \boldeta)$ and $\widehat{S}_{A} (t; \boldeta)$ are not smooth functions of $\boldeta$. In fact, they can be  very wiggly. As an illustration, we plot $\widehat{S}_{I} (t; \boldeta)$ and $\widehat{S}_{A} (t; \boldeta)$ as functions of $\eta_1$ in Figure 2 for a simple example with one covariate and the intercept of $\boldeta$ being set as 1. The black curves are for the estimates $\widehat{S}_{I} (t; \boldeta)$ and $\widehat{S}_{A} (t; \boldeta)$, which are given in the left and right panels of Figure 2, respectively. It can be clearly seen that the curves are very wiggly, and the direct maximization of them with respect to $\boldeta$ will be challenging and may lead to local maximizers. From our simulation studies conducted in Section~\ref{sec:simulation}, the estimated survival probability following the obtained optimal treatment regimes may have substantial biases.

\begin{figure}[ht]\label{value function}
\centering
\includegraphics[width=0.475\textwidth]{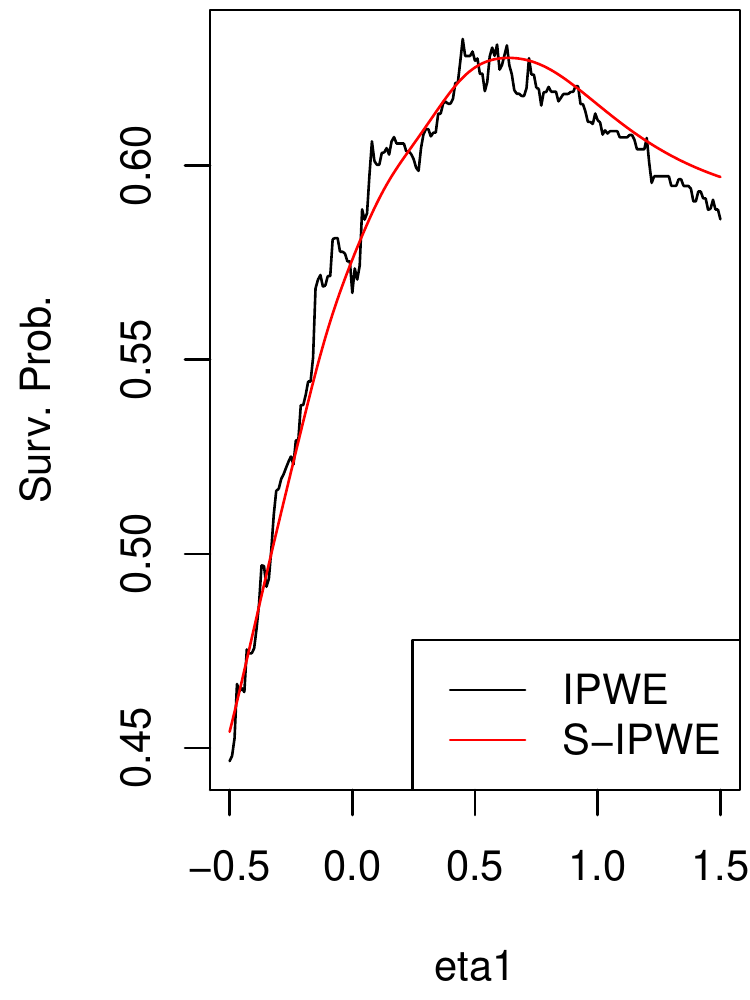}%
\includegraphics[width=0.475\textwidth]{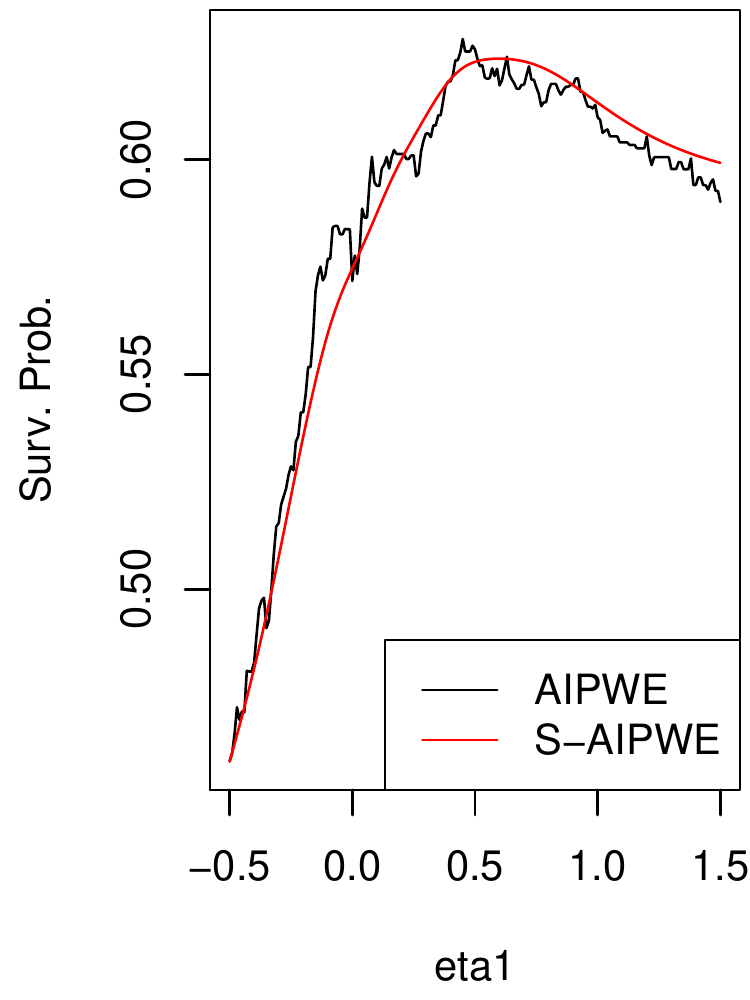}
\caption{Plots of original and smoothed value functions.}
\end{figure}

To reduce the biases of the estimates, we propose to smooth the estimates $\widehat{S}_{I} (t; \boldeta)$ and $\widehat{S}_{A} (t; \boldeta)$ using kernel smoothers. Specifically, we replace the $\geta(\boldX_i) = I\{ \boldeta^T \tilde{\boldX}_i \geq 0\}$ in $\widehat{S}_{I} (t; \boldeta)$ and $\widehat{S}_{A} (t; \boldeta)$  with $\tilde{g}_{\boldeta}(\boldX_i) = \Phi \left( \boldeta^T \tilde{\boldX}_i / h \right)$ to get the smoothed IPSWKME (S-IPSWKME) $\widetilde{S}_I(t; \boldeta)$ and the smoothed AIPSWKME (S-AIPSWKME) $\widetilde{S}_A(t; \boldeta)$, where $\Phi(s)$ is the cumulative distribution function for the standard normal distribution and $h$ is a bandwidth parameter  that goes to zero as $n$ goes to infinity. For the bandwidth selection, we set $h = c_0n^{1/3} \text{sd}(\boldeta^T\tilde{\boldX})$, where $c_0$ is a constant and $\text{sd}(\boldsymbol{v})$ is the sample standard deviation of $\boldsymbol{v}$. Such a bandwidth parameter has been widely used in nonparametric smoothing literature and will ensure that the original estimates and the smoothed estimates have the same asymptotic distributions. In our numerical studies, we tried different values for $c_0$ and found that $c_0 = 4^{1/3}$ generally gives good results for all scenarios. As an illustration, we plot in Figure 2 the smoothed estimates with the chosen bandwidth parameter for the same example in red curves. It can be seen that the smoothed curves well approximate the original curves and have unique maximizers around the true value $\eta_1 = 0.5$.
Let $\tilde{\boldeta}_I^{\text{opt}}$ and $\tilde{\boldeta}_A^{\text{opt}}$ denote the maximizers of  $\widetilde{S}_I(t; \boldeta)$ and $\widetilde{S}_A(t; \boldeta)$, respectively. Then the associated optimal treatment regimes are  $g(\boldx;\tilde{\boldeta}_I^{\text{opt}})$ and $g(\boldx;\tilde{\boldeta}_A^{\text{opt}})$.

\section{Estimation of Optimal Treatment Regime for Multiple Decision Time Points} \label{sec:methodology 2}

In this section, we extend our estimation methods to derive optimal dynamic treatment regimes  incorporating multiple decision time points. For the simplicity of presentation, we use the case with two decision time points as an illustration. Specifically, treatments can be given at the baseline and an interim time point $s$, $0 < s < t$.  For the $i\text{th}$ patient, let $\boldX_{0i} $ denote his or her $p_0$-dimensional vector of baseline covariates and  $A_{0i} \in \mathcal{A}_0 = \{0, 1\}$ denote the initial treatment received at the baseline. If this patient survives beyond $s$ and is not censored before $s$, let $\boldX_{1i}$ denote his or her $p_1$-dimensional vector of intermediate covariates collected by $s$ after assigning treatment $A_{0i}$ and $A_{1i} \in \mathcal{A}_1 = \{0,1\}$ denote the follow-up treatment given at $s$. Thus, the observed data are $\{(\boldX_{0i}, A_{0i}, \boldX_{1i}I\{\tilde{T}_i>u\}), A_{0i}I\{\tilde{T}_i > u\}, \tilde{T}_i, \delta_i), i=1,\ldots,n\}$.

As for single decision time point, we consider a class of linear dynamic treatment regimes for simplicity, i.e. $\mathcal{G}= \{ \boldsymbol{g}_{\boldeta}  = (g_0, g_1)\}$, where
\begin{align*}
g_0(\boldX_0; \boldeta_0) &= I\{\boldeta_{0}^T (1, \boldX_0^T) \geq 0\} ,\\
g_1(\boldX_0, \boldX_1; \boldeta_1) &= I\{ \boldeta_{1}^T (1, \boldX_0^T, g_0(\boldX_0; \boldeta_0), \boldX_1^T)) \geq 0\}, 
\end{align*}
and $\boldeta_0 \in \mathbb{R}^{p_0+1}, \boldeta_1 \in \mathbb{R}^{p_0 + p_1 +2}$. Here a patient following a treatment regime $\boldsymbol{g}_{\boldeta}$ implies that this patient is given treatment $g_0(\boldX_0; \boldeta_0) $ at baseline, and if he or she survives beyond $s$ and is not censored before $s$, this patient will be given treatment $g_1(\boldX_0, \boldX_1; \boldeta_1)$ at $s$. Note that for patients whose initial treatments coincide with those assigned by the regime $g_0(\boldX_0; \boldeta_0)$ and who die before $s$,  their treatment assignments are also consistent with the regime $\boldsymbol{g}_{\boldeta}$. However, for patients whose initial treatments coincide with those assigned by the regime $g_0(\boldX_0; \boldeta_0)$ but who are censored before $s$, it is not known whether their treatment assignments follow the regime $\boldsymbol{g}_{\boldeta}$.   Let $T^*(\boldsymbol{g}_{\boldeta}(\boldX_0,\boldX_1))$ denote the potential survival time for a patient if he or she were given treatment regime $\boldsymbol{g}_{\boldeta}(\boldX_0,\boldX_1)$. Here we are interested in finding the optimal dynamic treatment regime $\boldsymbol{g}_{\boldeta}^{\text{opt}} = (g_0(\boldX_0; \boldeta^{\text{opt}}_0),g_1(\boldX_0, \boldX_1; \boldeta^{\text{opt}}_1)) $ in  $\mathcal{G}$ that maximizes the $t$-year survival probability $S^{\ast(2)}(t; \boldeta)=E_{\boldX_0,\boldX_1}[P\{\Tstar(\geta(\boldX_0,\boldX_1)) > t | \boldX_0,\boldX_1\}]$. As commonly used in the causal inference literature for studying dynamic treatment regimes (e.g., Murphy, 2003), we make two assumptions: (i) SUTVA, i.e. a patient's observed outcome agrees with the corresponding potential outcome if his or her actually received treatments are consistent with the assigned treatments and (ii) sequential randomization assumption (SRA), i.e. the treatment assignment at current stage only depends on the past received treatments and observed covariates, but not the potential outcomes.  Under these two assumptions, the above defined $t$-year survival probability can be estimated from observed data.

Next, we propose a similar inverse propensity score weighted Kaplan-Meier estimator for the survival function $S^{\ast(2)}(u; \boldeta)$ given any treatment regime $\boldsymbol{g}_{\boldeta}$. However, the derivation of proper weights becomes more difficult since some patients may be censored before $s$ and whether their received treatments follow the regime $\boldsymbol{g}_{\boldeta}$ is unknown. To take this into account, we define the following new weight for patient $i$, $i=1,\ldots,n$:
\begin{align*}
\hatwetai^{(2)} = & \frac{I(\tilde{T}_i \le s )\times \delta_i }{\hat{S}_C(\tilde{T}_i)} \times \frac{I\{ A_{0i} = g_0(\boldX_{0i}; \boldeta_0)\}}{\hat{\pi}_{A_0}(\boldX_{0i})} \\
& + \frac{I(\tilde{T}_i > s)}{\hat{S}_C(s)} \times \frac{I\{A_{0i} = g_0(\boldX_{0i}; \boldeta_0), A_{1i} = g_1(\boldX_{0i}, g_0(\boldX_{0i};\boldeta_0), X_{1i}; \boldeta_1)\}}{\hat{\pi}_{A_0}(\boldX_{0i}) \times \hat{\pi}_{A_1}(\boldX_{0i}, \boldX_{1i})},
\end{align*}
where $\hat{\pi}_{A_0}(\boldX_{0i}) = \hat{\pi}_0(\boldX_{0i})A_{0i} + \{1-\hat{\pi}_0(\boldX_{0i})\}(1-A_{0i})$, $\hat{\pi}_{A_1}(\boldX_{0i},\boldX_{1i}) = \hat{\pi}_1(\boldX_{0i},\boldX_{1i})A_{1i} + \{1-\hat{\pi}_1(\boldX_{0i},\boldX_{1i})\}(1-A_{1i})$, and $\hat{\pi}_{0}(\boldX_{0i})$ and $\hat{\pi}_{1}(\boldX_{0i},\boldX_{1i}) $ are the estimates of the propensity scores $P(A_{0i}=1|\boldX_{0i})$ and $P(A_{1i}=1|\boldX_{0i},\boldX_{1i})$, respectively. In randomized studies, $\hat{\pi}_0$ and $\hat{\pi}_1$ are known by design, while in observational studies, they need to be estimated from day, say using logistic regression. Then the new IPSWKME for  $S^{\ast}(u; \boldeta)$  is given by
\begin{equation}\label{IPW-2}
\widehat{S}^{(2)}_I(u; \boldsymbol{\eta}) = \prod_{v \leq u} \left \{1- \frac{\sum_{i=1}^{n} \hatwetai^{(2)} d N_i(v)}{\sum_{i=1}^{n} \hatwetai^{(2)} Y_i(v)} \right \}.
\end{equation} 
Let $\hat{\boldsymbol{\eta}}_I^{\text{opt},(2)} = (\hat{\boldsymbol{\eta}}_{I,0}^{\text{opt},(2)},\hat{\boldsymbol{\eta}}_{I,1}^{\text{opt},(2)})= \arg\max_{||\boldsymbol{\eta}_0||=1,||\boldsymbol{\eta}_1||=1} \widehat{S}^{(2)}_I(t; \boldsymbol{\eta})$.  Then the estimated optimal dynamic treatment regime is given by $\hat{\boldsymbol{g}}_{\boldeta}^{\text{opt},(2)} = (g_0(\boldX_0; \hat{\boldeta}^{\text{opt},(2)}_{I,0}),g_1(\boldX_0, \boldX_1; \hat{\boldeta}^{\text{opt},(2)}_{I,1})) $.

To improve the finite sample performance of the IPSWKME, we also introduce kernel smoothing here. Specifically, we replace the indicator functions $g_0(\boldX_{0i}; \boldeta_0)$ and $g_1(\boldX_{0i}, \boldX_{1i}; \boldeta_1)$  in $\widehat{S}^{(2)}_I(u; \boldsymbol{\eta}) $ with $\Phi \left( \boldeta_{0}^T (1, \boldX_{0i}^T)/ h_0 \right)$ and 
$\Phi \left(  \boldeta_{1}^T (1, \boldX_0^T, g_0(\boldX_0; \boldeta_0), \boldX_1^T))/ h_1 \right)$, where the bandwidth parameters $h_0$ and $h_1$ are chosen similarly as before. Let $\widetilde{S}^{(2)}_I(u; \boldsymbol{\eta})$ denote the resulting smoothed IPSWKME and $\tilde{\boldsymbol{\eta}}_I^{\text{opt},(2)}$ denote the maximizer of $\widetilde{S}^{(2)}_I(t; \boldsymbol{\eta})$. To improve the robustness of IPSWKME, we can similarly derive the augmented IPSWKME based on a posited model for survival time, however, its formulation will be very complicated and is not pursued here. In addition, conceptually, the proposed IPSWKME can be generalized to accommodate more than two decision time points. However, when there are more treatment decision time points, the IPSWKME may become less reliable since fewer patients will follow a given dynamic treatment regime.

\section{Asymptotic Properties}\label{sec:asymptotic}
In this Section, we present the asymptotic properties of the proposed estimators which are summarized in Theorems 1 - 3. 

\begin{theorem}
Under conditions (A1)-(A6) in the Appendix, if the propensity score model~\eqref{eq:logistic} is correctly specified, for any regime $\geta$, we have, as $n \rightarrow \infty$,
\begin{enumerate}
\item[(i.)] $\widehat{S}_I(u; \boldeta) \rightarrow^{p} S^{\ast}(u; \boldeta)$ for any $0 <u \le t$; 
\item[(ii.)] $\sqrt{n}\{\widehat{S}_I(u; \boldeta) - S^{\ast}(u; \boldeta)\}$ converges weakly to a mean zero Gaussian process; 
\item[(iii.)] $\sqrt{n}\{\widehat{S}_I(t; \hat{\boldeta}_I^{\mathrm{opt}})-S^{\ast}(t; \boldeta^{\mathrm{opt}})\} \rightarrow^{d}N(0, \Sigma_I(t; \boldeta^{\mathrm{opt}}))$, where the expression of $\Sigma_I(t; \boldeta^{\mathrm{opt}})$ is given in the Appendix; 
\item[(iv.)] $\sqrt{n}\{\widehat{S}_I (t; \hat{\boldeta}_I^{\mathrm{opt}}) - \widetilde{S}_I (t; \tilde{\boldeta}_I^{\mathrm{opt}})\} = o_p(1)$.
\end{enumerate}
\end{theorem}

\begin{theorem}
Under condition (A1)-(A6) in the Appendix, if either the propensity score model~\eqref{eq:logistic} or the proportional hazard model~\eqref{eq:PH} is correctly specified, we have, as $n \rightarrow \infty$,
\begin{enumerate}
\item[(i.)] $\widehat{S}_A(u; \boldeta) \rightarrow^{p} S^{\ast}(u; \boldeta)$ for any $0 <u \le t$; 
\item[(ii.)] $\sqrt{n}\{\widehat{S}_A(u; \boldeta) - S^{\ast}(u; \boldeta)\}$ converges weakly to a mean zero Gaussian process; 
\item[(iii.)] $\sqrt{n}\{\widehat{S}_A(t; \hat{\boldeta}_A^{\mathrm{opt}})-S^{\ast}(t; \boldeta^{\mathrm{opt}})\} \rightarrow^{d}N(0, \Sigma_A(t; \boldeta^{\mathrm{opt}}))$, where the expression of $\Sigma_A(t; \boldeta^{\mathrm{opt}})$ is given in the Appendix; 
\item[(iv.)] $\sqrt{n}\{\widehat{S}_A(t; \hat{\boldeta}_A^{\mathrm{opt}}) - \widetilde{S}_A (t; \tilde{\boldeta}_A^{\mathrm{opt}})\} = o_p(1)$.
\end{enumerate}
\end{theorem}

\begin{theorem}
Under certain regularity conditions, if the two propensity score models $\pi_0(\cdot)$ and $\pi_1(\cdot)$ are correctly specified, for any regime $\geta$, we have, as $n \rightarrow \infty$,
\begin{enumerate}
\item[(i.)] $\widehat{S}^{(2)}_I(u; \boldeta) \rightarrow^{p} S^{\ast(2)}(u; \boldeta)$ for any $0 <u \le t$; 
\item[(ii.)] $\sqrt{n}\{\widehat{S}^{(2)}_I(u; \boldeta) - S^{\ast(2)}(u; \boldeta)\}$ converges weakly to a mean zero Gaussian process; 
\item[(iii.)] $\sqrt{n}\{\widehat{S}^{(2)}_I(t; \hat{\boldeta}_I^{\mathrm{opt},(2)})-S^{\ast}(t; \boldeta^{\mathrm{opt},(2)})\} \rightarrow^{d}N(0, \Sigma^{(2)}_I(t; \boldeta^{\mathrm{opt,(2)}}))$, where $\boldeta^{\mathrm{opt,(2)}} = (\boldeta_0^{\mathrm{opt}},\boldeta_1^{\mathrm{opt}})$;
\item[(iv.)] $\sqrt{n}\{\widehat{S}^{(2)}_I (t; \hat{\boldeta}_I^{\mathrm{opt,(2)}}) - \widetilde{S}^{(2)}_I (t; \tilde{\boldeta}_I^{\mathrm{opt,(2)}})\} = o_p(1)$.
\end{enumerate}
\end{theorem}

Here the asymptotic variance $\Sigma_I(t; \boldeta^{\mathrm{opt}})$, $\Sigma_A(t; \boldeta^{\mathrm{opt}})$ and $\Sigma^{(2)}_I(t; \boldeta^{\mathrm{opt,(2)}})$ can be consistently estimated from observed data using the usual plug-in method. The proofs of Theorems 1-3 are given in the Appendix.

\section{Simulation Studies}\label{sec:simulation}

In this Section, we examine the finite sample performance of the proposed estimators by simulations. We first consider scenarios with a single treatment decision time point at the baseline. For each patient, the baseline covariates $X_1$ and $X_2$ are independently and uniformly distributed on $(-2,2)$. Given the covariates $X_1$ and $X_2$, the binary treatment indicator $A$ is generated from the logistic model $\text{logit} \{\pi(X_1, X_2)\} = X_1 - 0.5 X_2$. The survival time $T$ is generated from a linear transformation model \citep{Linear-Transformation-Model}, $h(T) = -0.5X_1+A(X_1-X_2) + \varepsilon$, where $h(s) = \log(e^s-1)-2$ is an increasing function and the error term $\varepsilon$ follows some known distribution, taking either the extreme value distribution or the logistic distribution, which corresponds to a proportional hazards and proportional odds model, respectively. The covariate-independent censoring time $C$ is uniformly distributed on $(0, C_0)$, where $C_0$ is chosen to achieve the censoring rate of $15\%$ and $40\%$. 
It is obvious the optimal treatment regime for maximizing $t$-year survival probability is $g_{\boldeta}^{\text{opt}}(X_1, X_2) = I \{ X_1-X_2 \geq 0 \}$ for any $t$. Here, we search the optimal treatment regime in the class of regimes given by $\mathcal{G} = \{ \geta: \geta(X_1, X_2)= I\{ \eta_0+ \eta_1 X_1 + \eta_2 X_2 \geq 0\}, \boldeta \in \mathbb{R}^3 \}$,  which contains the true optimal treatment regime as a special case. For easy comparion, we impose the restriction $\boldeta^T \boldeta = 1$ and thus we have $\boldeta^{\text{opt}}=(0, 0.707, -0.707)$.

To implement our proposed estimators, we need to posit a model for the propensity scores. Here, we consider both a correctly specified model: $\text{logit} \{\pi(X_1, X_2)\} = \theta_0 + \theta_1 X_1 + \theta_2 X_2$ and a misspecified model: $\text{logit} \{\pi_A(X_1, X_2)\} = \theta_0$. For the augmented estimators, we need to posit a model for the survival time $T$. Here, we always use the proportional hazard model $\lambda(t | X_1, X_2) = \lambda_0(t) \exp \{\beta_{11} X_1 + \beta_{12}X_2 + A(\beta_{20} + \beta_{21} X_1+ \beta_{22} X_2)\}$. Note that when $\varepsilon$ follows the extreme value distribution, the posited survival model is correctly specified. On the other hand, when $\varepsilon$ follows the logistic distribution, this model is misspecified. We compared the performance of the IPSWKME ($\hat{S}_I$) and AIPSWKME ($\hat{S}_A$), as well as their smoothed versions: S-IPSWKME ($\tilde{S}_I$) and S-AIPSWKME ($\tilde{S}_A$), under different combinations of the assumed propensity score (PS) model, error term distribution, censoring rate, sample size ($n = 250$ or 500) and time point of interest ($t = 1$ or 2). For each scenario, we run 1000 replications and use the genetic algorithm to do the optimization, which is implemented by the \texttt{R} function \texttt{genoud} within the package \texttt{rgenoud} \citep{Genetic-Algorithm}.

To save the presentation space, we only report the simulation results for the scenarios with $n = 250$ and $t = 2$, which are given in Tables 1 and 2 for the extreme value error  and logistic error distributions, respectively. Results for other scenarios are very similar and omitted here. In the tables, we report the mean of estimated $\boldeta$, the mean of estimated $t$-year survival probability following the estimated optimal treatment regime, namely the estimated optimal $t$-year survival probability (denoted by $\hat{S}(\hat{\boldeta}^{\text{opt}}))$, the mean of estimated standard error of $\hat{S}(\hat{\boldeta}^{\text{opt}})$ using the plug-in method based on the asymptotic variances established in Theorems 1-2 (denoted by SE), the empirical coverage probability of 95\% confidence interval for the $t$-year survival probability following the true optimal treatment regime $S(\boldeta^{\text{opt}})$ (denoted by CP), the mean of simulated true $t$-year survival probability following the estimated optimal treatment regime (denoted by $S(\hat{\boldeta}^{\text{opt}})$), and the mean of misclassification rate by comparing the true and estimated optimal treatment regimes (denoted by MR). The numbers given in parenthesis are the standard deviation of the corresponding estimates. Here, $S(\boldeta^{\text{opt}})$ and $S(\hat{\boldeta}^{\text{opt}})$ are computed using simulated survival times following the given treatment regime based on a large random sample of $5 \times 10^6$ patients. We have $S(\boldeta^{\text{opt}}) = 0.605$ for the extreme value error distribution and $S(\boldeta^{\text{opt}}) = 0.672$ for the logistic distribution. In addition, the misclassification rate for one simulation is calculated as the proportion of patients that the true and estimated optimal treatment regimes do not match. 


From the results, we make the following observations. First, when the PS model is correctly specified, all the estimators of $\boldeta^{\text{opt}}$ have relatively small biases, in particular, the  mean of $\hat{\eta}_0^{\text{opt}}$ is close to zero while the mean ratio of $\hat{\eta}_1^{\text{opt}}$ to $\hat{\eta}_2^{\text{opt}}$ is very close to negative one. The means of simulated true $t$-year survival probability following the estimated optimal treatment regimes, i.e. $S(\hat{\boldeta}^{\text{opt}})$, are all close to the true values.  In addition, the  estimates of $\boldeta^{\text{opt}}$ based on the AIPSWKME and S-AIPSWKME of $t$-year survival probability generally have smaller standard deviation than those based on IPSWKME and S-IPSWKME. Second, the unsmoothed IPSWKME and AIPSWKME of the optimal $t$-year survival probability have relatively large biases mainly due to the very wiggly estimates of $t$-year survival probability as illustrated in Figure 2 and as a consequence, the associated coverage probability of 95\% confidence interval is much lower than the nominal level. Third, the smoothed S-IPSWKME and S-AIPSWKME of the optimal $t$-year survival probability greatly reduce the biases and thus give the proper coverage probability. In addition, the unsmoothed and smoothed estimators of the optimal $t$-year survival probability have nearly the same standard deviation. 
Fourth, when the PS model is misspecified, the IPSWKME and S-IPSWKME generally have relatively large biases as expected, while the AIPSWKME and S-AIPSWKME greatly reduce the biases and give much smaller MR. In particular, when the posited survival model is correctly specified under the extreme value error distribution, the S-AIPSWKME gives proper coverage probability. On the other hand, when the posited survival model is misspecified under the logistic error distribution, although the S-AIPSWKME is not consistent in general, it still gives small biases with reasonable coverage probability. Lastly, the performance of our proposed estimators improve as the censoring rate decreases and sample size increases.


Next, we consider scenarios with two treatment decision time points, one at the baseline and the other at $s=1$. The initial treatment assignment $A_0$ and the follow-up treatment assignment $A_1$, if applicable, are generated independently from a Bernoulli distribution with success probability of 0.5. A single baseline covariate is generated from a uniform distribution on $(0, 4)$. To generate the survival time $T$, we first generate a time $T_1$ given $A_0$ and $X_0$ from an exponential distribution with the rate function $\lambda_1(A_0, X_0)$. The censoring time $C$ is generated from a uniform distribution on $(0, C_0)$. If a patient is neither dead nor censored at time $s=1$ (i.e. $\min (T_1, C) > 1$), we generate a single intermediate covariate $X_1$ for this patient by $X_1 = 0.5X_0 - 0.4 (A_0-0.5) + e$, where $e$ is uniformly distributed on $(0,2)$. Then we generate another time $T_2$ given $A_0$, $A_1$, $X_0$ and $X_1$ from an exponential distribution with the rate function $\lambda_2(A_0, A_1, X_0, X_1)$. The survival time $T$ of interest is defined as $T = T_1$ if $T_1 \le 1$ and $T = 1 + T_2$ otherwise. The observed survival time is $\tilde{T} = \min(T,C)$ with the censoring indicator $\delta = I(T \le C)$. Here the constant $C_0$ is chosen to achieve the censoring rate of $15\%$ and $40\%$. We consider three scenarios for the rate functions $\lambda_1$ and $\lambda_2$: (i) $\lambda_1(A_0,X_0) = 0.5 \exp \{1.75(A_0-0.5)(X_0 - 2)\}$ and $\lambda_2(A_0, A_1, X_0, X_1) = 0.3 \exp\{2.5(A_1-0.4)(X_1-2) - A_0(X_1-2)\}$; (ii) $\lambda_1(A_0,X_0) = 0.1 \exp \{2(A_0-0.5)(X_0 - 2)\}$ and $\lambda_2(A_0, A_1, X_0, X_1) = 0.2 \exp\{3(A_1-0.4)(X_1-2) - 3(A_0-0.5)(X_0-2)\}$; (iii) $\lambda_1(A_0,X_0) = 0.2 \exp \{1.5 (A_0-0.3) (X_0 - 3)\}$ and $\lambda_2(A_0, A_1, X_0, X_1) = 0.3 \exp\{2 (A_1-0.5) (X_1-2) +0.5 (A_0-0.7)(X_0-1)\}$. 

For the above three scenarios, it is easy to see that the true optimal treatment regime for maximizing $t$-year survival probability ($t > 1$) at time $s=1$ is given by $g_1^{\text{opt}} = I(2-X_1 > 0)$. However, the true optimal treatment regime $g_0^{\text{opt}}$ at time $s=0$ is a very complicated nonlinear function of $X_0$, which can be derived using backward induction as done in Q-learning. In our implementation, for computation simplicity, we search the optimal dynamic treatment regime in a class of linear decision rules, specifically, $\mathcal{G}_{\boldeta} = \{g_0(X_0) = I\{\eta_1 + \eta_2 X_0 >0\}, g_1(X_1) = I\{\eta_3 + \eta_4 X_1 >0\}, ||(\eta_1, \eta_2)||=1, ||(\eta_3, \eta_4)||=1\}$. It is clear that the true optimal treatment regime at $s=1$ is contained in the class but 
the true optimal treatment regime at $s=0$ is not. For scenarios (i) and (iii), we take $t= 3$, while for (ii) we take $t = 6$. Instead of finding the true optimal treatment regime at $s=0$, we use simulation method to find the best treatment regime at $s=0$ in the class $\mathcal{G}_{\boldeta} $ to maximize $t$-year survival probability. To be specific, we first generate $X_0$, and for a given $(\eta_1, \eta_2)$, we set $A_0$ by the regime $g_0(X_0)$. Then, we generate $X_1$ given $A_0$ and $X_0$ the same way as in our design, and set $A_1$ by the optimal regime $g_1^{\text{opt}}$. Finally, we generate $T_1$ and $T_2$, and define $T$ the same way as before. Based on the generated $T$'s for a large random sample of $5 \times 10^6$ patients, we compute the associated empirical $t$-year survival probability. We find $(\eta_1^{\text{opt}}, \eta_2^{\text{opt}})$ to maximize the empirical $t$-year survival probability, which gives the best treatment regime $g_0^{\text{opt}}$ in the class $\mathcal{G}_{\boldeta} $. Here we use grid search method to find $(\eta_1^{\text{opt}}, \eta_2^{\text{opt}})$. Since $||(\eta_1^{\text{opt}}, \eta_2^{\text{opt}})|| = 1$, we only need to do grid search for $\eta_1$. We have $(\eta_1^{\text{opt}}, \eta_2^{\text{opt}}) = (0.890,-0.456)$ and $S(3; \etaopt)=0.567$ for scenario 1, $(\eta_1^{\text{opt}}, \eta_2^{\text{opt}}) = (-0.891,0.454)$ and $S(6; \etaopt)=0.624$ for scenario 2, and $(\eta_1^{\text{opt}}, \eta_2^{\text{opt}}) = (0.908,-0.419)$ and $S(3; \etaopt)=0.702$ for scenario 3. Here $\boldeta^{\text{opt}} = (\eta_1^{\text{opt}}, \eta_2^{\text{opt}}, \eta_3^{\text{opt}}, \eta_4^{\text{opt}})$ and $S(t;\etaopt)$ is the $t$-year survival probability following the optimal dynamic treatment regime $\etaopt$. Note that $(\eta_3^{\text{opt}}, \eta_4^{\text{opt}}) = (0.894,-0.447)$ after normalization for all three scenarios.  

We compare the unsmoothed and smoothed estimators. For both estimators, the propensity score models $\pi_0$ and $\pi_1$ are assumed known as for randomized clinical trials. Simulation results are summarized in Table 3. From the results, we observe: (i) both unsmoothed and smoothed estimation methods give nearly unbiased estimators of $\etaopt$, and the $t$-year survival probability following the estimated optimal treatment regime (denoted by $S(\hat{\boldeta}^{\text{opt}})$ in the table) is very close to the $t$-year survival probability following the true optimal treatment regime $\etaopt$; (ii) the mean of estimated standard error (SE) of  $\hat{S}(\hat{\boldeta}^{\text{opt}})$ based on the established theory is close to the standard deviation of the estimates given in the parenthesis;  (iii) The unsmoothed estimator for the $t$-year survival probability following the estimated optimal treatment regime (denoted by $\hat{S}(\hat{\boldeta}^{\text{opt}})$) has relatively large bias and the associated coverage probability (CP) is below the nominal level; and (iv) the smoothed estimator for the $t$-year survival probability following the estimated optimal treatment regime has largely reduced bias and thus lead to proper coverage probability.

\section{A Data Example}\label{sec:application}
We illustrate the proposed methods with the data from the AIDS Clinical Trials Group Study 175 \citep{Real-Data}. This is a randomized clinical trial and patients were randomized to four treatment groups with equal probability: zidovudine (ZDV) monotherapy, ZDV plus didanosine (ddI), ZDV plus zalcitabine (zal), and ddI
monotherapy. A primary endpoint of interest is the time to having a larger than 50\% decline in the CD4 count, or progressing to AIDS, or death, whichever comes first. From treatment-specific Kaplan-Meier curves, it can be clearly seen that treatments ZDV+ddI, ZDV+zal and ddI only are uniformly better than treatment ZDV only in terms of survival. In addition, treatments ZDV+ddI and ZDV+zal are overall the two best treatments giving the highest survival probabilities especially after day 400.
For simplicity,  we only consider two treatment options in our analysis , specifically, $A=1$ for zidovudine+ddI and $A=0$ for zidovudine+zal, which involves 1046 patients. For each patient, there are 12 baseline clinical covariates. From historical studies (e.g., Geng et al., 2014), it is found that Karnofsky score (Karnof), baseline CD4 count (CD40), and age (Age) are three important risk predictors and may have interaction effects with treatments. In our analysis, we only include these three covariates in constructing treatment regimes. Our goal is to find the optimal treatment regime TO from the class of linear regimes defined by $\mathcal{G}=\{ g_{\boldeta}=I(\eta_0 + \eta_1 \text{Karnof} + \eta_2 \text{CD40}+ \eta_3 \text{Age} \geq 0) :  \boldeta \in \mathbb{R}^4 \}$ to maximize $t$-year survival probability. To simplify notation, we define $X_1$ as Karnof, $X_2$  as CD40 and $X_3$ as Age. Since the data comes from a randomized study, we use a constant model for the propensity score and estimate this constant from data. For the augmented estimation,  we posit the proportional hazard model as given in (5).  We consider $t=400$, 600, 800 and 1000. We only compute the S-IPSWKME and S-AIPSWKME, since they have better numerical performance than their nonsmoothed counterparts based on our simulation studies. 

The estimated optimal treatment regimes and the associated $t$-year survival probabilities are presented in Table~\ref{tab:realdata}. The numbers given in the columns of  Intercept, Karnof, CD40 and Age are the parameter estimates $\tilde{\boldeta}^{\text{opt}}$ defining the optimal treatment regimes, and $\tilde{S}(t;\tilde{\boldeta}^{\text{opt}})$ is the estimated $t$-year survival probability following the estimated optimal treatment regime. We make the following observations: (i) the estimated optimal treatment regime at earlier time may be different from that at later time. For example, comparing the obtained optimal treatment regimes at $t = 600$ and $t = 800$, the S-IPSWKME assigns a set of 355 patients to treatment 0 and another set of 585 patients to treatment 1 at both time points. However, it assigns a set of 51 patients to treatment 0 at day 600 but to treatment 1 at day 800. On the other hand, it assigns another set of 55 patients to treatment 1 at day 600 but to treatment 0 at day 800. For the S-AIPSWKME, the findings are similar. (ii) The S-IPSWKME and S-AIPSWKME may give very different parameter estimates $\tilde{\boldeta}^{\text{opt}}$. However, the corresponding optimal treatment regimes may be similar. Using the results at day 600 as an example, among the 1046 patients, there are only 68 patients whose assigned treatments are different by the estimated optimal treatment regimes based on S-IPSWKME and S-AIPSWKME. In addition, the estimated $t$-year survival probabilities following the estimated optimal treatment regimes are nearly the same based on S-IPSWKME and S-AIPSWKME. 

Next, we compare the estimated optimal regimes with the simple regimes that assign everyone to the same treatment. Specifically, we construct the 95\% confidence intervals for the difference between the estimated $t$-year survival probabilities under the estimated optimal treatment regimes and the simple regimes using two methods: one is the Wald-type confidence interval based on the derived asymptotic normal distribution and the other is the bootstrap confidence interval based on 500 runs. The results are given in Table~\ref{tab:compare}. From the results we observe that (i) the Wald-type confidence interval and bootstrap confidence interval are very similar; (ii) the bootstrap confidence intervals all stay above 0 when comparing the estimated $t$-year survival probabilities under the estimated optimal treatment regimes and the simple regimes  for all the considered time points, indicating that the estimated optimal treatment regimes significantly improves $t$-year survival probabilities comparing with simple regimes; (iii) Some Wald-type confidence interval based on normal approximation stays above 0 and others contain 0. However, for those that contain 0, zero is very close to the left end of the intervals. Therefore, similar conclusions can be made here as for the bootstrap confidence intervals, although they are a little less significant.

\section{Discussion}\label{sec:discussion}
In this paper, we propose various Kaplan-Meier type estimators for the survival function of patients following a given (dynamic) treatment regime. We further introduce kernel smoothing for the proposed estimators to improve their numerical performance. Then, the optimal (dynamic) treatment regime is searched within a class of pre-specified treatment regimes to maximize the associated $t$-year survival probability. Current work only considers the case when there are two treatment options at each decision time point. However, the proposed method can be generalized to incorporate multiple treatment options at each decision time point by defining a treatment regime using multiple indexes instead of a single indicator function $g_{\boldeta}(\boldX) =I\{\boldeta^T \tilde{\boldX} \geq 0\}$.  In addition, current methods find the optimal (dynamic) treatment regime to maximize $t$-year survival probability, which can also be generalized to maximize other clinical outcomes of interest. Specifically, using the IPSWKME, $\widehat{S}_I(\cdot;\boldeta)$, as an illustration,  we can find the optimal treatment regime to maximize $f\{\widehat{S}_I(\cdot;\boldeta)\}$, where $f$ is a specified function of interest. For example, if we take $f\{\widehat{S}_I(\cdot;\boldeta)\} = \int_0^L\widehat{S}_I(u;\boldeta)du $, which corresponds to the restricted mean survival time under a given treatment regime. On the other hand, if we take $f\{\widehat{S}_I(\cdot;\boldeta)\} = \sup\{u: \widehat{S}_I(u;\boldeta) \ge 0.5\} $ , which corresponds to the median survival time under a given treatment regime. These are interesting topics that need further investigation.


\appendix

\setcounter{equation}{0}
\renewcommand{\theequation}{A.\arabic{equation}}

\section{Proof of Theorems}

To establish the asymptotic results given in Theorems 1-2, we need to assume some regularity conditions. Recall that a working logistic model (3) is assumed for the propensity scores with parameters $\boldtheta$ for the IPSWKME and a working proportional hazards model (5) is further assumed for the survival time $T$ for the AIPSWKME with parameters $\boldbeta$ and $\Lambda_0$.  Let $\nuAi = (\boldX_i^T, A_i, A_i \boldX_i^T)^T$ and $\nuetai = (\boldX_i^T, \geta(\boldX_i), \geta(\boldX_i) \boldX_i^T)^T$. Define 
\begin{align*}
&K_1^I(\boldX, A, \tilde{T}, \delta; \boldeta) = \int_{0}^{t} \frac{(2A-1)dN(u)}{\pi^{\ast}E\{w_{\boldeta}^{\ast} Y(u) \}}, \\
&K_2^I(\boldX, A, \tilde{T}, \delta; \boldeta) =  \int_{0}^{t} \frac{(2A-1)Y(u) E[\{ (2A-1) \geta(\boldX)  + (1-A) \}dN(u)]}{[\pi^{\ast}E\{w_{\boldeta}^{\ast} Y(u) \}]^2}, 
\end{align*}
where $w_{\boldeta}^{\ast} = [A\geta(\boldX) + (1-A) \{1-\geta(\boldX)\}]/\pi^{\ast}$ and $\pi^{\ast} = \pi(\boldX; \boldtheta^{\ast}) A + \{1-\pi(\boldX; \boldtheta^{\ast})\} (1-A)$.
In addition, define
\begin{align*}
&K_1^A(\boldX, A, \tilde{T}, \delta; \boldeta) = \int_0^{t} \frac{J_1^A(u) - J_0^A(u) }{E\left[ \{L_1^A(u) - L_0^A(u)\} {g}_{\boldeta}(\boldX) + L_0^A(u) \right]}, \\
&K_2^A(\boldX, A, \tilde{T}, \delta; \boldeta) =\int_0^{t} \frac{\{L_1^A(u) -L_0^A(u)\} E\left[ \{J_1^A(u) -J_0^A(u)\} \geta(\boldX) + J_0^A(u)\right]}{ \left( E\left[ \{L_1^A(u) - L_0^A(u)\} \ {g}_{\boldeta}(\boldX) + L_0^A(u) \right] \right)^2},
\end{align*}
where $J_k^A(u) = \frac{1-k - (-1)^kA}{\pi^{\ast}} d N(u) + e_k\left(1-\frac{1-k - (-1)^kA}{\pi^{\ast}}\right) \exp\left\{-\Lambda_0^{\ast}(u)e_k\right\} S_C(u) d \Lambda_0^{\ast}(u)$,  $L_k^A(u) = \frac{1-k - (-1)^kA}{\pi^{\ast}} Y(u) + \left(1-\frac{1-k - (-1)^kA}{\pi^{\ast}}\right) \exp\left\{-\Lambda_0^{\ast}(u)e_k\right\} S_C(u)$, $e_k = \exp\left\{ {\boldbeta^{\ast}}^T (\boldX^T, k, k\boldX^T)^T\right\}$, $k =0, 1$. We assume the following conditions.
\begin{enumerate}
\item[A1.] The covariates $\boldX$ are bounded.
\item[A2.] The propensity score $\pi(\boldX)$ is bounded away from $0$ and $1$ for all possible values of $\boldX$.
\item[A3.] The equation $E\left[ \left\{ A - \frac{\exp (\boldtheta^T \tilde{\boldX})}{1+ \exp(\boldtheta^T \tilde{\boldX})} \right\} \tilde{\boldX} \right] = 0$ has a unique solution $\boldtheta^{\ast}$.
\item[A4.] The equation 
\begin{equation*}
E \left( \int_0^{\tau} \left[\nuAi - \frac{E\left\{Y_i(s) \exp(\boldbeta^T \nuAi) \nuAi \right\}}{E \left\{Y_i(s) \exp(\boldbeta^T \nuAi )\right\}} \right] \times d N_i(s)  \right) =0.
\end{equation*}
has a unique  solution $\boldbeta^{\ast}$, where $\tau > t$ is a pre-specified time point satisfying $P(\tilde{T}_i \geq \tau) > 0$. Let $\Lambda_0^*(u) = E[\int_0^udN_i(s)/E \{Y_i(s) \exp({\boldbeta^*}^T \nuAi )\}]$ and it satisfies $\Lambda_0^*(\tau) < \infty$. 
\item[A5.] $\text{sup}_{||\boldeta||=1} E[\{K_j^I(\boldX, A, \tilde{T}, \delta; \boldeta)\}^2] < \infty$ and $\text{sup}_{||\boldeta||=1} E[\{K_j^A(\boldX, A, \tilde{T}, \delta; \boldeta)\}^2] < \infty$,  $j =1, 2$. 
\item[A6.] $nh \rightarrow \infty$ and $nh^4 \rightarrow 0$ as $n \rightarrow \infty$.
\end{enumerate}

Under assumed regularity conditions A1 - A4, we have the following asymptotic representations:
\begin{equation*}
\sqrt{n}(\hat{\boldtheta} - \boldtheta^{\ast}) =  \frac{1}{\sqrt{n}} \sum_{i=1}^{n} \phi_{1 i} + o_p(1), \quad \sqrt{n}(\hat{\boldbeta} - \boldbeta^{\ast}) = \frac{1}{\sqrt{n}} \sum_{i=1}^{n} \phi_{2 i} + o_p(1),
\end{equation*}
\begin{equation*}
\sqrt{n}\{\hat{\Lambda}_0(u) - \Lambda_0^{\ast}(u)\} = \frac{1}{\sqrt{n}} \sum_{i=1}^{n} \phi_{3 i} (u) + o_p(1), \quad \sqrt{n}\{\hat{S}_C(u) - S_C(u)\} = \frac{1}{\sqrt{n}} \sum_{i=1}^{n} \phi_{4 i}(u) + o_p(1),
\end{equation*}
where $ \phi_{1 i}$'s and  $\phi_{2 i}$'s are independently and identically distributed mean-zero vectors, and $\phi_{3 i} (u)$ and $\phi_{4 i} (u)$ are independent mean-zero processes.

\subsection{Proof of Theorem 1}
For any given regime $g_{\boldeta}$, we first derive the asymptotic properties for the corresponding inverse propensity score weighted (IPSW) Nelson-Aalen estimator. Specifically,
\begin{equation}\label{eq:IPWE}
\widehat{\Lambda}_I(u; \boldeta)  \equiv \widehat{\Lambda}_I(u; \boldeta,\hat{\boldtheta})= \int_0^{u} \frac{\sum_{i=1}^{n} \hatwetai d N_i(s)}{ \sum_{i=1}^{n} \hatwetai Y_i(s)}.
\end{equation}
It is easy to show that $\widehat{S}_I(u; \boldeta)$ and $\exp\{-\widehat{\Lambda}_I(u; \boldeta) \}$ are asymptotically equivalent for any given $\boldeta$.
Therefore, the asymptotic properties of $\widehat{S}_I(u; \boldeta)$ easily follows those of $\widehat{\Lambda}_I(u; \boldeta)$.

When the propensity score model is correctly specified, we have that $\boldtheta^{\ast} = \boldtheta$ and $\wetai^{\ast}=\wetai$. Then $n^{-1} \sum_{i=1}^{n} \hatwetai Y_i(s)
\rightarrow_p E\{\wetai Y_i(s)\} = E[Y^*\{g_{\boldeta}(X);s\}]$ uniformly for $s \in [0, \tau]$ as $n \rightarrow \infty$. Similarly, we have $n^{-1} \sum_{i=1}^{n} \hatwetai dN_i(s)
\rightarrow_p E\{\wetai dN_i(s)\} = E[dN^*\{g_{\boldeta}(X);s\}]$ uniformly for $s \in [0, \tau]$ as $n \rightarrow \infty$. Therefore, 
\begin{eqnarray*}
\widehat{\Lambda}_I(u; \boldeta)  \rightarrow_p & &\int_{0}^{u} \frac{E[dN^*\{g_{\boldeta}(X);s\}]}{E[Y^*\{g_{\boldeta}(X);s\}]} = \int_{0}^{u} \frac{S_C(s)dP[T^*\{g_{\boldeta}(X)\}\le s]}{S_C(s)P[T^*\{g_{\boldeta}(X)\}\ge s]}\\
& &  = -\log\{S^*(u;\boldeta)\} \equiv \Lambda^*(u;\boldeta),
\end{eqnarray*}
which establish the consistency given in (i) of Theorem 1.

Next, we derive the asymptotic distribution of $\widehat{\Lambda}_I(u; \boldeta)$. By applying the first-order Taylor expansion of  $\widehat{\Lambda}_I(u; \boldeta)$ with respect to parameter $\boldtheta$, we have
$$
\sqrt{n}\{\widehat{\Lambda}_I(u; \boldeta) - \Lambda^*(u; \boldeta)\}
= \sqrt{n}\{\widehat{\Lambda}_I(u; \boldeta,\boldtheta) - \Lambda^*(u; \boldeta)\} +  D_1(u)^T\sqrt{n}(\hat{\boldtheta} - \boldtheta)  + o_p(1),
$$
where $D_1(u) = \lim_{n \rightarrow \infty}\partial \widehat{\Lambda}_I(u; \boldeta,\boldtheta)/\partial \boldtheta$. In addition,
\begin{eqnarray*}
&& \sqrt{n}\{\widehat{\Lambda}_I(u; \boldeta,\boldtheta) - \Lambda^*(u; \boldeta)\}  = \sqrt{n}\int_0^{u} \frac{\sum_{i=1}^{n} \wetai \{d N_i(s) - Y_i(s)d\Lambda^*(s; \boldeta)\}}{ \sum_{i=1}^{n} \wetai Y_i(s)}\\
&=&n^{-1/2}\sum_{i=1}^n\int_0^{u} \frac{\wetai [d N_i^*\{g_{\boldeta}(X);s\} - Y_i^*\{g_{\boldeta}(X);s\}d\Lambda^*(s; \boldeta)]}{E[Y^*\{g_{\boldeta}(X);s\}]}+o_p(1)\\
&=& n^{-1/2}\sum_{i=1}^n\int_0^{u} \frac{\wetai d M_i^*\{g_{\boldeta}(X);s\}}{E[Y^*\{g_{\boldeta}(X);s\}]}+o_p(1),
\end{eqnarray*}
where $M_i^*\{g_{\boldeta}(X);s\} = N_i^*\{g_{\boldeta}(X);s\} - \int_0^s Y_i^*\{g_{\boldeta}(X);v\}d\Lambda^*(v; \boldeta)$ is a mean-zero martingale process. Therefore,
\begin{eqnarray*}
\sqrt{n}\{\widehat{\Lambda}_I(u; \boldeta) - \Lambda^*(u; \boldeta)\}
&=& n^{-1/2}\sum_{i=1}^n\left(\int_0^{u} \frac{\wetai d M_i^*\{g_{\boldeta}(X);s\}}{E[Y^*\{g_{\boldeta}(X);s\}]} + D_1(u)^T\phi_{1i}\right)+o_p(1)\\
&\equiv& n^{-1/2}\sum_{i=1}^n \zeta_i(u;\boldeta) + o_p(1),
\end{eqnarray*}
where $\zeta_i(u;\boldeta)$'s are independent mean-zero processes. 
By delta method, we have
$\sqrt{n} \{ \widehat{S}_I(u; \boldeta) - S^*(u; \boldeta) \} = -S^*(u; \boldeta)n^{-1/2}\sum_{i=1}^n \zeta_i(u;\boldeta) + o_p(1)$, 
which converges weakly to a mean-zero Gaussian process by applying the empirical process theory. This proves (ii) of Theorem 1.

Since $\hatetaopt_I$ is the maximizer of $\widehat{S}_I(t; \boldeta)$ and $\etaopt$ is the maximizer of $S^*(t; \boldeta)$, following the similar arguments in \citet{Baqun-Zhang}, we have
$$
\sqrt{n}\{\widehat{S}_I(t; \hat{\boldeta}_I^{\text{opt}}) - S^*(t; \boldeta^{\text{opt}}) \} -\sqrt{n} \{ \widehat{S}_I(t; {\boldeta}^{\text{opt}}) - S^*(t; \boldeta^{\text{opt}}) \} = o_p(1).
$$
It follows that $\sqrt{n}\{\widehat{S}_I(t; \hat{\boldeta}_I^{\mathrm{opt}})-S^{\ast}(t; \boldeta^{\mathrm{opt}})\} \rightarrow^{d}N(0, \Sigma_I(t; \boldeta^{\mathrm{opt}}))$, where 
$\Sigma_I(t; \boldeta^{\mathrm{opt}}) = \{S^*(u; \boldeta^{\mathrm{opt}})\}^2E\{\zeta_i^2(u;\boldeta^{\mathrm{opt}})\}$.
This proves (iii) of Theorem 1.

Finally, we show that $\widehat{S}_I(t; \hat{\boldeta}_I^{\mathrm{opt}})$ and $\widetilde{S}_I(t; \tilde{\boldeta}_I^{\mathrm{opt}})$ are asymptotically equivalent. 
 For any given $\boldeta$, we have
\begin{align}
& \sqrt{n}\left\{\widetilde{\Lambda}_I(t; \boldeta) - \widehat{\Lambda}_I(t; \boldeta) \right\} \notag \\
&=  \sqrt{n} \times \frac{1}{n} \sum_{i=1}^{n} \left\{ \Phi \left( \frac{\boldeta^T \boldX_i}{h} \right) -I\left( \boldeta^T \boldX_i \geq 0 \right)  \right\} \times  K_1^I(\boldX_i, A_i, \tilde{T}_i, \delta; \boldeta) \label{eq:t1} \\
&+  \sqrt{n} \times \frac{1}{n} \sum_{i=1}^{n} \left\{ \Phi \left( \frac{\boldeta^T \boldX_i}{h} \right) -I\left( \boldeta^T \boldX_i \geq 0 \right)  \right\} \times  K_2^I(\boldX_i, A_i, \tilde{T}_i, \delta; \boldeta) \label{eq:t2}\\
&+o_p(1), \notag
\end{align}
For simplicity, define $\boldsymbol{q} = (\boldX_i, A_i, \tilde{T}_i, \delta)$ and $r^{\boldeta}=\boldeta^T \boldX$ . Following the similar arguments in \citet{Smoothing}, we have
\begin{equation*}
| \eqref{eq:t1} | \leq M \sqrt{n} \text{ sup}_{||\boldeta||=1} \left|\int_{\boldsymbol{q}} \int_{r^{\boldeta}}\left\{ \Phi\left( \frac{r^{\boldeta}}{h} \right) - I(r^{\boldeta} \geq 0)\right\} K_1^I(\boldsymbol{q}; \boldeta) d \hat{F}(r^{\boldeta} | \boldsymbol{q}; \boldeta) d \hat{G}(\boldsymbol{q}; \boldeta)\right|, \label{eq:B1}
\end{equation*}
where $M$ is a finite constant, $\hat{G}(\boldsymbol{q}; \boldeta)$ and $\hat{F}(r^{\boldeta} | \boldsymbol{q}; \boldeta)$ are the marginal empirical cumulative distribution functions for $\boldsymbol{q}$ and the conditional empirical cumulative distribution function for $r^{\boldeta}$, respectively. For simplicity, we omit the superscript $\boldeta$ in $r^{\boldeta}$, the condition $\boldeta$ in $K_1^I(\boldsymbol{q}; \boldeta)$, $\hat{F}(r|\boldsymbol{q}; \boldeta)$ and $\hat{G}(\boldsymbol{q}; \boldeta)$. Thus, the equation~\eqref{eq:t1} is bounded by $M \sqrt{n}\sup_{||\boldeta||=1} | \Upsilon|$, where
$$
\Upsilon = \int_{\boldsymbol{q}} \int_r \left\{ \Phi\left( \frac{r}{h} \right) - I(r \geq 0)\right\} K_1^I(\boldsymbol{q}) d \hat{F}(r | \boldsymbol{q}) d \hat{G}(\boldsymbol{q}).
$$
Write $\Upsilon = \Upsilon_1 + \Upsilon_2$, where
\begin{align*}
\Upsilon_1 &=  \int_{\boldsymbol{q}} \int_r \left\{ \Phi\left( \frac{r}{h} \right) - I(r \geq 0)\right\} K_1^I(\boldsymbol{q}) \left\{ d \hat{F}(r | \boldsymbol{q}) - d F(r|\boldsymbol{q}) \right\} d \hat{G}(\boldsymbol{q}) \\
\Upsilon_2 &=  \int_{\boldsymbol{q}} \int_r \left\{ \Phi\left( \frac{r}{h} \right) - I(r \geq 0)\right\} K_1^I(\boldsymbol{q}) d F(r | \boldsymbol{q}) d \hat{G}(\boldsymbol{q})
\end{align*}
with $F(r|\boldsymbol{q})=\lim_{n \to +\infty} \hat{F}(r|\boldsymbol{q})$. By variable transformation  $z=r/h$ and integration by parts, we have
\begin{equation}\label{eq:Upsilon1}
\Upsilon_1 = \int_{\boldsymbol{q}} \int_z K_1^I(\boldsymbol{q}) \varphi(z) \left\{ \left[ \hat{F}(zh|\boldsymbol{q}) - F(zh|\boldsymbol{q})\right] - \left[ \hat{F}(0|\boldsymbol{q}) - F(0|\boldsymbol{q}))\right] \right\} d z d\hat{G}(\boldsymbol{q}),
\end{equation}
where $\varphi(z)$ is the probability density function of standard normal distribution. Under regularity condition A5, we apply the results on oscillations of empirical process \citep{Empirical-Process} to equation~\eqref{eq:Upsilon1} and have
$$
\sqrt{n} |\Upsilon_1  | = O_p\left(\sqrt{h \log n \log \left(\frac{1}{h \log n} \right)}\right).
$$
In addition, by similar arguments and applying second order Taylor expansion of $\Upsilon_2$ with respect to $h$ around 0, we have
$$
\Upsilon_2 = -\frac{ h^2}{2} \int_{\boldsymbol{q}} \int_z K_1^I(\boldsymbol{q}) \varphi(z) f'(zh^{\ast}|\boldsymbol{q}) z^2 dz d\hat{G}(\boldsymbol{q}),
$$
where $f'(u|\boldsymbol{q}) = \partial^2 F(u|\boldsymbol{q}) / \partial u^2$ and $h^{\ast}$ lies between $h$ and $0$. Thus, we have
$\sqrt{n}|\Upsilon_2|= O_p(\sqrt{n}h^2)$.
Combine the above results, we have
$$
|\eqref{eq:t1}| \leq \sqrt{n}|\Upsilon_1|+\sqrt{n}|\Upsilon_2|  = O_p\left( \sqrt{h \log n \log \left(\frac{1}{h \log n} \right)}+  \sqrt{n}h^2 \right).
$$
By condition A6, we have $\sup_{||\boldeta||=1}|\eqref{eq:t1}| = o_p(1)$.
Similarly, we have $\sup_{||\boldeta||=1}|\eqref{eq:t2}| = o_p(1)$. Therefore, we have
$\sqrt{n}\{\widetilde{\Lambda}_I(t; \boldeta) -\widehat{\Lambda}_I(t; \boldeta)\} = o_p(1)$ uniformly in $\boldeta$, which implies 
$\sqrt{n}\{\widetilde{S}_I(t; \boldeta) -\widehat{S}_I(t; \boldeta)\} = o_p(1)$ uniformly in $\boldeta$. In addition, it is easy to show that $\sqrt{n}\{\widetilde{S}_I(t; \tilde{\boldeta}_I^{\text{opt}}) -\widetilde{S}_I(t; \boldeta^{\text{opt}})\} = o_p(1)$ and $\sqrt{n}\{\widehat{S}_I(t; \hat{\boldeta}_I^{\text{opt}}) -\widehat{S}_I(t; \boldeta^{\text{opt}})\} = o_p(1)$. It follows that $\sqrt{n}\{\widetilde{S}_I(t; \tilde{\boldeta}_I^{\text{opt}}) -\widehat{S}_I(t; \hat{\boldeta}_I^{\text{opt}})\} = o_p(1)$, which proves (iv) of Theorem 1.

\subsection{Proof of Theorem 2}

For any given regime $\geta$, we similarly introduce the augmented IPSW Nelson-Aalen estimator
\begin{equation}\label{eq:AIPWENA}
\widehat{\Lambda}_A (u; \boldeta)  = \int_{0}^{u}
\frac{\sum_{i=1}^n \hatwetai d N_i(s) + (1-\hatwetai) \hat{S}_T(s | \geta(\boldX_i), \boldX_i) \hat{S}_C(s) d\hat{\Lambda}_T (s | \geta(\boldX_i), \boldX_i)}
{\sum_{i=1}^n \hatwetai Y_i(s) + (1-\hatwetai) \hat{S}_T(s | \geta(\boldX_i), \boldX_i) \hat{S}_C(s)}.
\end{equation}
We will show that $\widehat{\Lambda}_A (u; \boldeta)$ is consistent when either the propensity score model is correctly specified or the survival model for $T$ is correctly specified, i.e. having the doubly robustness property. First, assume that the propensity score model is correctly specified. Then, we have $\boldtheta^{\ast} = \boldtheta$ and $\wetai^{\ast} = \wetai$. In addition, the denominator of equation~\eqref{eq:AIPWENA} converges in probability to 
$E\{\wetai Y_i(s)\}  + E \left[(1-\wetai)  \exp \{-\Lambda_0^{\ast}(s) \exp ({\boldbeta^{\ast}}^T \nu_{\boldeta i}) \}  S_C(s)\right]$ uniformly for $s \in [0,\tau]$. Note that the second term is zero since $E(\wetai|\boldX_i) = 0$. Similarly, the numerator of equation~\eqref{eq:AIPWENA} converges in probability to 
$$
E\{ \wetai d N_i(u)\}+ E \left[(1-\wetai) \exp \{ -\Lambda_0^{\ast}(u) \exp ({\boldbeta^{\ast}}^T \nu_{\boldeta i})\} S_C(u) \exp({\boldbeta^{\ast}}^T \nu_{\boldeta i}) d \Lambda_0^{\ast}(u) \right ] 
$$ 
uniformly for $s \in [0,\tau]$, where the second term is also zero. The proof of consistency then follows that for the IPSW Nelson-Aalen estimator.

On the other hand, when the survival model for $T$ is correctly specified, we have $\boldbeta^{\ast} = \boldbeta$ and $\Lambda_0^{\ast}(s) = \Lambda_0(s)$. We can show that
the denominator of equation~\eqref{eq:AIPWENA} converges in probability to 
$$
E\left[  \exp \{-\Lambda_0(s) \exp (\boldbeta^T \nu_{\boldeta i})\} S_C(s) \right] +E\left(\wetai^{\ast} [Y_i(s) - \exp \{-\Lambda_0(s) \exp (\boldbeta^T \nu_{\boldeta i}) \} S_C(s) ] \right)
$$
uniformly for $s \in [0,\tau]$, where the first term equals to $S^*(s; \boldeta) S_C(s)$ and the second term is zero since $E[Y_i(s) - \exp \{-\Lambda_0(s) \exp (\boldbeta^T \nu_{\boldeta i}) \} S_C(s) |A_i,\boldX_i] = 0$.
In addition, the numerator of equation~\eqref{eq:AIPWENA} converges in probability to 
\begin{align*}
&E\left[ \exp \{ -\Lambda_0(s) \exp ({\boldbeta}^T \nu_{\boldeta i})\} S_C(s) \exp({\boldbeta}^T \nu_{\boldeta i}) d \Lambda_0(s) \right]  \\
+ &E\left(\wetai^{\ast}  [d N_i(u) - \exp \{ -\Lambda_0(s) \exp ({\boldbeta}^T \nu_{g i})\} S_C(s) \exp({\boldbeta}^T \nu_{\boldeta i})  d \Lambda_0(u) ]\right) 
\end{align*}
uniformly for $s \in [0,\tau]$, where the first term equals to $-S_C(s)dS^*(s; \boldeta) $ and the second term is zero since $E [d N_i(u) - \exp \{ -\Lambda_0(s) \exp ({\boldbeta}^T \nu_{g i})\} S_C(s) \exp({\boldbeta}^T \nu_{\boldeta i})  d \Lambda_0(u)|A_i,\boldX_i ] = 0$.
Therefore, the remaining proof follows that for the IPSW Nelson-Aalen estimator.

Next, we derive the asymptotic distribution for $\widehat{S}_A(u; \boldeta)$, assuming that either the propensity score model or the survival model for $T$ is correctly specified. Note that $\widehat{\Lambda}_A (u; \boldeta) = \widehat{\Lambda}_A (u; \boldeta,\hat{\boldtheta},\hat{\boldbeta},\hat{\Lambda}_0,\hat{S}_C)$. By Taylor expansion of $\widehat{\Lambda}_A (u; \boldeta,\hat{\boldtheta},\hat{\boldbeta},\hat{\Lambda}_0,\hat{S}_C)$ with respect to the estimators $\hat{\boldtheta}$, $\hat{\boldbeta}$, $\hat{\Lambda}_0$ and $\hat{S}_C$ around their population values, we have 
$$
\sqrt{n}\{\widehat{\Lambda}_A(u; \boldeta) - \Lambda^*(u; \boldeta)\} = \sqrt{n}\{\widehat{\Lambda}_A(u; \boldeta,\boldtheta^*,\boldbeta^*,\Lambda^*_0,S_C) - \Lambda^*(u; \boldeta)\} + n^{-1/2}\sum_{i=1}^n \psi_{2i}(u;\boldeta) + o_p(1),
$$
where $\psi_2(u;\boldeta)$'s are independent mean-zero processes due to the asymptotic expansions of the estimators $\hat{\boldtheta}$, $\hat{\boldbeta}$, $\hat{\Lambda}_0$ and $\hat{S}_C$ .
By simple algebra, we have
\begin{equation*}
\sqrt{n}\{\widehat{\Lambda}_A(u; \boldeta,\boldtheta^*,\boldbeta^*,\Lambda^*_0,S_C) - \Lambda^*(u; \boldeta)\}
= n^{-1/2}\sum_{i=1}^{n} \int_{0}^{u}  \frac{d h_i(s) }{ E[Y^*\{g_{\boldeta}(X);s\}]} + o_p(1),
\end{equation*}
where
\begin{align*}
dh_i(s) =& \wetai^{\ast} \{d N_i(s) -  Y_i(s) d\Lambda^*(s; \boldeta)\}\\
&+ (1-\wetai^{\ast}) S^*_T(s | \geta(\boldX_i), \boldX_i) S_C(s)d \{\Lambda^*_T (s | \geta(\boldX_i), \boldX_i) -\Lambda^*(s; \boldeta)\}.
\end{align*}
Note that the first term in $dh_i(s)$ equals to $\wetai^{\ast}d M_i^*\{g_{\boldeta}(X);s\}$ and the second term is zero if the propensity score model is correctly specified. If the survival model for $T$ is correctly specified, we have $E\{\Lambda^*_T (s | \geta(\boldX_i), \boldX_i) \} - \Lambda^*(s; \boldeta) = 0$. Define $\psi_{1i} (u;\boldeta)= \int_{0}^{u}  \frac{d h_i(s) }{ E[Y^*\{g_{\boldeta}(X);s\}]} $. Then, $\psi_{1i} (u;\boldeta)$'s are independent mean-zero processes. Let $\psi_i(u;\boldeta) = \psi_{1i} (u;\boldeta) + \psi_{2i} (u;\boldeta)$.
We have $
\sqrt{n}\{\widehat{\Lambda}_A(u; \boldeta) - \Lambda^*(u; \boldeta)\}  = n^{-1/2}\sum_{i=1}^n \psi_{i}(u;\boldeta) + o_p(1)$, which converges weakly to a mean-zero Gaussian process. By Delta method, $
\sqrt{n}\{\widehat{S}_A(u; \boldeta) -S^*(u; \boldeta)\} $ also converges weakly to a mean-zero Gaussian process. 

Following the proof for Theorem 1, we have 
$$
\sqrt{n}\{\widehat{S}_A(t; \hat{\boldeta}_A^{\text{opt}}) - S^*(t; \boldeta^{\text{opt}}) \} -\sqrt{n} \{ \widehat{S}_A(t; {\boldeta}^{\text{opt}}) - S^*(t; \boldeta^{\text{opt}}) \} = o_p(1).
$$
It follows that $\sqrt{n}\{\widehat{S}_A(t; \hat{\boldeta}_A^{\mathrm{opt}})-S^{\ast}(t; \boldeta^{\mathrm{opt}})\} \rightarrow^{d}N(0, \Sigma_A(t; \boldeta^{\mathrm{opt}}))$, where 
$\Sigma_A(t; \boldeta^{\mathrm{opt}}) = \{S^*(u; \boldeta^{\mathrm{opt}})\}^2E\{\psi_i^2(u;\boldeta^{\mathrm{opt}})\}$.

Finally, for any given $\boldeta$, we have
\begin{align}
& \sqrt{n}\left\{\widetilde{\Lambda}_A(t; \boldeta) - \widehat{\Lambda}_A(t; \boldeta) \right\} \notag \\
&=  \sqrt{n} \times \frac{1}{n} \sum_{i=1}^{n} \left\{ \Phi \left( \frac{\boldeta^T \boldX_i}{h} \right) -I\left( \boldeta^T \boldX_i \geq 0 \right)  \right\} \times  K_1^A(\boldX_i, A_i, \tilde{T}_i, \delta; \boldeta) \label{eq:t3} \\
&+  \sqrt{n} \times \frac{1}{n} \sum_{i=1}^{n} \left\{ \Phi \left( \frac{\boldeta^T \boldX_i}{h} \right) -I\left( \boldeta^T \boldX_i \geq 0 \right)  \right\} \times  K_2^A(\boldX_i, A_i, \tilde{T}_i, \delta; \boldeta) \label{eq:t4}\\
&+o_p(1). \notag
\end{align}
Under conditions A5 and A6, following the similar arguments in the proof for (iv) of Theorem 1, (\ref{eq:t3}) and (\ref{eq:t4}) can be bounded uniformly in $\boldeta$. Therefore, $\sqrt{n}\{\widetilde{S}_A(t; \boldeta) -\widehat{S}_A(t; \boldeta)\} = o_p(1)$ uniformly in  $\boldeta$. Since $\sqrt{n}\{\widetilde{S}_A(t; \tilde{\boldeta}_A^{\text{opt}}) -\widetilde{S}_A(t; \boldeta^{\text{opt}})\} = o_p(1)$ and $\sqrt{n}\{\widehat{S}_A(t; \hat{\boldeta}_A^{\text{opt}}) -\widehat{S}_A(t; \boldeta^{\text{opt}})\} = o_p(1)$, it follows that $\sqrt{n}\{\widetilde{S}_A(t; \tilde{\boldeta}_A^{\text{opt}}) -\widehat{S}_A(t; \hat{\boldeta}_A^{\text{opt}})\} = o_p(1)$.

\subsection{Proof of Theorem 3}

To establish the asymptotic results given in Theorem 3, the regularity conditions A1-A3 and A5-A6 need to be modified accordingly to incorporate the two-stage treatment regimes, and condition A4 is not needed. However, the proof of Theorem 3 can follow similar steps as for the proof of Theorem 1, and is omitted here.

\bibliographystyle{apalike}
\bibliography{project2}

\newpage

\begin{landscape}

\begin{table}[ht]
\centering
\caption{Simulation results for the extreme value error distribution with $n= 250$ and  $t = 2$.} \smallskip
\begin{tabular}{ccrccccccc}
  \hline
\hline
  & PS & \ccell{$\hat{\eta}_0$} & $\hat{\eta}_1$ & $\hat{\eta}_2$ & $\hat{S}(\hat{\boldeta}^{\text{opt}})$ & SE & CP & $S(\hat{\boldeta}^{\text{opt}})$ & MR \\
  \hline \multicolumn{10}{c}{Censoring Rate = 15\%} \\$\hat{S}_I$ & T & $0.008\;(0.302)$ & $0.631\;(0.191)$ & $-0.666\;(0.179)$ & $0.645\;(0.037)$ & $0.040$ & $0.839$ & $0.590\;(0.016)$ & $0.118\;(0.064)$ \\
  $\tilde{S}_I$ & T & $-0.005\;(0.262)$ & $0.653\;(0.179)$ & $-0.666\;(0.171)$ & $0.612\;(0.036)$ & $0.040$ & $0.968$ & $0.593\;(0.014)$ & $0.107\;(0.057)$ \\
  $\hat{S}_A$ & T & $0.006\;(0.285)$ & $0.639\;(0.172)$ & $-0.675\;(0.161)$ & $0.639\;(0.037)$ & $0.041$ & $0.882$ & $0.592\;(0.014)$ & $0.109\;(0.059)$ \\
  $\tilde{S}_A$ & T & $-0.002\;(0.260)$ & $0.654\;(0.175)$ & $-0.670\;(0.160)$ & $0.610\;(0.036)$ & $0.041$ & $0.970$ & $0.593\;(0.013)$ & $0.104\;(0.056)$ \\
  $\hat{S}_I$ & F & $-0.026\;(0.414)$ & $0.413\;(0.321)$ & $-0.702\;(0.249)$ & $0.666\;(0.036)$ & $0.039$ & $0.657$ & $0.566\;(0.038)$ & $0.190\;(0.099)$ \\
  $\tilde{S}_I$ & F & $-0.051\;(0.402)$ & $0.427\;(0.284)$ & $-0.714\;(0.252)$ & $0.643\;(0.035)$ & $0.039$ & $0.844$ & $0.569\;(0.034)$ & $0.184\;(0.090)$ \\
  $\hat{S}_A$ & F & $-0.013\;(0.277)$ & $0.661\;(0.152)$ & $-0.662\;(0.160)$ & $0.635\;(0.038)$ & $0.041$ & $0.889$ & $0.593\;(0.011)$ & $0.106\;(0.055)$ \\
  $\tilde{S}_A$ & F & $0.001\;(0.315)$ & $0.616\;(0.183)$ & $-0.669\;(0.200)$ & $0.612\;(0.037)$ & $0.042$ & $0.966$ & $0.589\;(0.015)$ & $0.126\;(0.062)$ \\
   \hline \multicolumn{10}{c}{Censoring Rate = 40\%} \\$\hat{S}_I$ & T & $0.004\;(0.317)$ & $0.615\;(0.215)$ & $-0.659\;(0.202)$ & $0.650\;(0.041)$ & $0.044$ & $0.848$ & $0.587\;(0.019)$ & $0.128\;(0.069)$ \\
  $\tilde{S}_I$ & T & $-0.002\;(0.286)$ & $0.637\;(0.202)$ & $-0.660\;(0.192)$ & $0.613\;(0.040)$ & $0.045$ & $0.958$ & $0.590\;(0.017)$ & $0.118\;(0.064)$ \\
  $\hat{S}_A$ & T & $0.003\;(0.305)$ & $0.621\;(0.204)$ & $-0.664\;(0.199)$ & $0.645\;(0.041)$ & $0.046$ & $0.892$ & $0.589\;(0.019)$ & $0.124\;(0.067)$ \\
  $\tilde{S}_A$ & T & $0.002\;(0.290)$ & $0.642\;(0.196)$ & $-0.656\;(0.188)$ & $0.612\;(0.040)$ & $0.046$ & $0.966$ & $0.590\;(0.017)$ & $0.118\;(0.064)$ \\
  $\hat{S}_I$ & F & $-0.002\;(0.439)$ & $0.394\;(0.344)$ & $-0.677\;(0.275)$ & $0.671\;(0.040)$ & $0.043$ & $0.678$ & $0.561\;(0.043)$ & $0.204\;(0.106)$ \\
  $\tilde{S}_I$ & F & $-0.024\;(0.432)$ & $0.404\;(0.310)$ & $-0.694\;(0.271)$ & $0.645\;(0.039)$ & $0.043$ & $0.867$ & $0.564\;(0.038)$ & $0.199\;(0.094)$ \\
  $\hat{S}_A$ & F & $-0.005\;(0.302)$ & $0.652\;(0.168)$ & $-0.650\;(0.183)$ & $0.641\;(0.042)$ & $0.046$ & $0.894$ & $0.591\;(0.014)$ & $0.116\;(0.061)$ \\
  $\tilde{S}_A$ & F & $0.011\;(0.339)$ & $0.606\;(0.204)$ & $-0.655\;(0.217)$ & $0.615\;(0.041)$ & $0.046$ & $0.961$ & $0.586\;(0.018)$ & $0.138\;(0.067)$ \\
   \hline
\end{tabular}
\label{tab:E1t2N250}
$^\dagger$ PS, the propensity score model. Here T means the correctly specified PS model while F means the misspecified PS model. Recall that $S(\boldeta^{\text{opt}}) = 0.605$.
\end{table}

\begin{table}[ht]
\centering
\caption{Simulation results for the logistic error distribution with $n= 250$ and  $t = 2$.} \smallskip
\begin{tabular}{ccrccccccc}
  \hline
\hline
  & PS & \ccell{$\hat{\eta}_0$} & $\hat{\eta}_1$ & $\hat{\eta}_2$ & $\hat{S}(\hat{\boldeta}^{\text{opt}})$  & SE & CP & $S(\hat{\boldeta}^{\text{opt}})$ & MR \\
  \hline \multicolumn{10}{c}{Censoring Rate = 15\%} \\$\hat{S}_I$ & T & $0.013\;(0.374)$ & $0.559\;(0.277)$ & $-0.641\;(0.246)$ & $0.716\;(0.034)$ & $0.038$ & $0.790$ & $0.652\;(0.023)$ & $0.156\;(0.092)$ \\
  $\tilde{S}_I$ & T & $-0.002\;(0.340)$ & $0.593\;(0.259)$ & $-0.641\;(0.235)$ & $0.685\;(0.034)$ & $0.039$ & $0.955$ & $0.655\;(0.020)$ & $0.145\;(0.081)$ \\
  $\hat{S}_A$ & T & $0.008\;(0.360)$ & $0.576\;(0.257)$ & $-0.645\;(0.235)$ & $0.713\;(0.034)$ & $0.040$ & $0.833$ & $0.654\;(0.020)$ & $0.149\;(0.084)$ \\
  $\tilde{S}_A$ & T & $-0.009\;(0.343)$ & $0.592\;(0.256)$ & $-0.642\;(0.233)$ & $0.684\;(0.034)$ & $0.040$ & $0.964$ & $0.655\;(0.020)$ & $0.144\;(0.082)$ \\
  $\hat{S}_I$ & F & $0.033\;(0.462)$ & $0.342\;(0.388)$ & $-0.662\;(0.284)$ & $0.729\;(0.033)$ & $0.037$ & $0.649$ & $0.632\;(0.039)$ & $0.223\;(0.119)$ \\
  $\tilde{S}_I$ & F & $-0.002\;(0.460)$ & $0.376\;(0.350)$ & $-0.666\;(0.285)$ & $0.707\;(0.033)$ & $0.037$ & $0.846$ & $0.636\;(0.034)$ & $0.216\;(0.107)$ \\
  $\hat{S}_A$ & F & $-0.019\;(0.336)$ & $0.627\;(0.203)$ & $-0.638\;(0.213)$ & $0.723\;(0.036)$ & $0.040$ & $0.757$ & $0.658\;(0.013)$ & $0.134\;(0.068)$ \\
  $\tilde{S}_A$ & F & $-0.022\;(0.353)$ & $0.594\;(0.224)$ & $-0.646\;(0.234)$ & $0.698\;(0.035)$ & $0.040$ & $0.920$ & $0.656\;(0.015)$ & $0.146\;(0.070)$ \\
   \hline \multicolumn{10}{c}{Censoring Rate = 40\%} \\$\hat{S}_I$ & T & $0.013\;(0.385)$ & $0.548\;(0.293)$ & $-0.630\;(0.261)$ & $0.721\;(0.036)$ & $0.041$ & $0.784$ & $0.650\;(0.026)$ & $0.165\;(0.095)$ \\
  $\tilde{S}_I$ & T & $-0.007\;(0.361)$ & $0.581\;(0.273)$ & $-0.626\;(0.256)$ & $0.687\;(0.036)$ & $0.041$ & $0.948$ & $0.652\;(0.022)$ & $0.155\;(0.087)$ \\
  $\hat{S}_A$ & T & $0.008\;(0.379)$ & $0.559\;(0.277)$ & $-0.632\;(0.261)$ & $0.718\;(0.036)$ & $0.043$ & $0.814$ & $0.651\;(0.023)$ & $0.160\;(0.090)$ \\
  $\tilde{S}_A$ & T & $-0.018\;(0.360)$ & $0.578\;(0.271)$ & $-0.634\;(0.247)$ & $0.687\;(0.036)$ & $0.043$ & $0.961$ & $0.653\;(0.022)$ & $0.153\;(0.086)$ \\
  $\hat{S}_I$ & F & $0.048\;(0.472)$ & $0.329\;(0.411)$ & $-0.635\;(0.307)$ & $0.733\;(0.035)$ & $0.039$ & $0.658$ & $0.628\;(0.042)$ & $0.236\;(0.125)$ \\
  $\tilde{S}_I$ & F & $0.020\;(0.481)$ & $0.358\;(0.367)$ & $-0.638\;(0.314)$ & $0.709\;(0.035)$ & $0.040$ & $0.842$ & $0.631\;(0.038)$ & $0.229\;(0.113)$ \\
  $\hat{S}_A$ & F & $-0.005\;(0.349)$ & $0.620\;(0.207)$ & $-0.636\;(0.217)$ & $0.722\;(0.038)$ & $0.043$ & $0.788$ & $0.657\;(0.015)$ & $0.138\;(0.071)$ \\
  $\tilde{S}_A$ & F & $-0.010\;(0.376)$ & $0.581\;(0.239)$ & $-0.634\;(0.250)$ & $0.696\;(0.038)$ & $0.043$ & $0.932$ & $0.653\;(0.016)$ & $0.156\;(0.074)$ \\
   \hline
\end{tabular}
\label{tab:E2t2N250}
$^\dagger$ PS, the propensity score model. Here T means the correctly specified PS model while F means the misspecified PS model. Recall that $S(\boldeta^{\text{opt}}) = 0.672$.
\end{table}

\begin{table}[ht]
\centering
\caption{Simulation results for estimating optimal dynamic treatment regimes.} \smallskip
\begin{tabular}{ccccccccccc}
  \hline \hline
$C\%$ & S & \ccell{$\hat{\eta}_1^{\text{opt}}$} &  \ccell{$\hat{\eta}_2^{\text{opt}}$} &  \ccell{$\hat{\eta}_3^{\text{opt}}$} &  \ccell{$\hat{\eta}_4^{\text{opt}}$} &  \ccell{$\hat{S}(\hat{\boldeta}^{\text{opt}})$} &  \ccell{SE} &  \ccell{CP} &  \ccell{$S(\hat{\boldeta}^{\text{opt}})$} &  \ccell{MR} \\ 
  \hline
\multicolumn{11}{c}{Senario 1: $\etaopt=(0.890, -0.456, 0.894, -0.447); S(3; \etaopt)=0.567$} \\
15 & F & $0.882\;(0.035)$ & $-0.466\;(0.062)$ & $0.893\;(0.016)$ & $-0.449\;(0.032)$ & $0.591\;(0.028)$ & $0.030$ & $0.885$ & $0.559\;(0.008)$ & $0.105\;(0.054)$ \\ 
& T & $0.884\;(0.028)$ & $-0.463\;(0.052)$ & $0.894\;(0.013)$ & $-0.448\;(0.026)$ & $0.570\;(0.028)$ & $0.030$ & $0.955$ & $0.561\;(0.006)$ & $0.088\;(0.048)$ \\ 
40 & F & $0.880\;(0.041)$ & $-0.469\;(0.071)$ & $0.890\;(0.022)$ & $-0.453\;(0.041)$ & $0.600\;(0.036)$ & $0.037$ & $0.841$ & $0.556\;(0.011)$ & $0.124\;(0.061)$ \\ 
& T  & $0.883\;(0.03)$ & $-0.463\;(0.061)$ & $0.892\;(0.018)$ & $-0.450\;(0.035)$ & $0.574\;(0.035)$ & $0.038$ & $0.955$ & $0.558\;(0.009)$ & $0.108\;(0.056)$ \\ 
   \hline
\multicolumn{11}{c}{Senario 2: $\etaopt=(-0.891, 0.454, 0.894, -0.447); S(6; \etaopt)=0.624$} \\
15 & F & $-0.888\;(0.025)$ & $0.456\;(0.044)$ & $0.891\;(0.018)$ & $-0.451\;(0.034)$ & $0.645\;(0.025)$ & $0.027$ & $0.890$ & $0.616\;(0.008)$ & $0.097\;(0.051)$ \\ 
& T & $-0.889\;(0.018)$ & $0.456\;(0.034)$ & $0.893\;(0.014)$ & $-0.450\;(0.028)$ & $0.624\;(0.024)$ & $0.027$ & $0.967$ & $0.618\;(0.005)$ & $0.079\;(0.042)$ \\ 
 40 & F & $-0.886\;(0.028)$ & $0.460\;(0.051)$ & $0.891\;(0.020)$ & $-0.453\;(0.037)$ & $0.650\;(0.027)$ & $0.029$ & $0.857$ & $0.614\;(0.009)$ & $0.108\;(0.054)$ \\ 
 & T & $-0.888\;(0.022)$ & $0.457\;(0.040)$ & $0.892\;(0.016)$ & $-0.450\;(0.032)$ & $0.626\;(0.027)$ & $0.030$ & $0.972$ & $0.617\;(0.007)$ & $0.091\;(0.048)$ \\ 
\hline
\multicolumn{11}{c}{Senario 3: $\etaopt=(0.908, -0.419, 0.894, -0.447); S(3; \etaopt)=0.702$} \\
15 & F & $0.898\;(0.037)$ & $-0.433\;(0.068)$ & $0.892\;(0.020)$ & $-0.450\;(0.038)$ & $0.728\;(0.026)$ & $0.027$ & $0.829$ & $0.693\;(0.009)$ & $0.132\;(0.067)$ \\ 
 & T & $0.900\;(0.031)$ & $-0.430\;(0.060)$ & $0.893\;(0.016)$ & $-0.448\;(0.031)$ & $0.707\;(0.026)$ & $0.027$ & $0.952$ & $0.695\;(0.007)$ & $0.115\;(0.060)$ \\ 
 40 & F & $0.897\;(0.040)$ & $-0.435\;(0.074)$ & $0.891\;(0.022)$ & $-0.452\;(0.042)$ & $0.732\;(0.028)$ & $0.029$ & $0.808$ & $0.691\;(0.011)$ & $0.140\;(0.074)$ \\ 
 & T & $0.899\;(0.035)$ & $-0.431\;(0.065)$ & $0.893\;(0.018)$ & $-0.449\;(0.036)$ & $0.709\;(0.028)$ & $0.030$ & $0.951$ & $0.693\;(0.008)$ & $0.125\;(0.065)$ \\ 
   \hline
\end{tabular}
$^\dagger$$C\%$ denotes the censoring rate; $S$ indicates whether the smoothing technique is applied (T) or not (F). 
\end{table}

\begin{table}[ht]
\centering
\caption{Estimation results for the AIDS data.}\smallskip
\begin{tabular}{ccrrrrl}
  \hline\hline
 $t$ &  \ccell{Method}  & \ccell{Intercept} & \ccell{Karnof} & \ccell{CD40} & \ccell{Age} & \ccell{$\tilde{S}(t;\tilde{\boldeta}^{\text{opt}})$} \\
  \hline
$400$ & I & -0.143 & -0.355 & 0.025 & 0.924 & $0.965\;(0.008)$ \\
   & A  & -0.660 & -0.265 & 0.020 & 0.703 & $0.965\;(0.008)$ \\
  $600$ & I & 0.908 & -0.147 & 0.002 & 0.391 & $0.923\;(0.012)$ \\
   & A& 0.998 & -0.026 & -0.000 & 0.050 & $0.923\;(0.012)$ \\
  $800$ & I & 0.815 & -0.154 & -0.011 & 0.558 & $0.887\;(0.014)$ \\
   & A & 0.882 & -0.127 & -0.009 & 0.453 & $0.886\;(0.014)$ \\
  $1000$ & I& 0.067 & -0.192 & -0.035 & 0.978 & $0.824\;(0.017)$ \\
   & A& -0.619 & -0.140 & -0.029 & 0.772 & $0.823\;(0.018)$ \\ 
   \hline 
\end{tabular}

$^\dagger$I denotes the IPSWKME and A denotes the AIPSWKME; the numbers in the parenthesis are the estimated standard errors. 
\label{tab:realdata}
\end{table}

\end{landscape}

\begin{table}[ht]
\centering
\caption{Confidence intervals for comparing estimated optimal treatment regimes and simple regimes.} \smallskip
\begin{tabular}{cccccc}
  \hline\hline
  &    & \multicolumn{2}{c}{Norm CI} & \multicolumn{2}{c}{Boot CI} \\
  \cline{3-4}  \cline{5-6}
 $t$ &  Method  & trt 1 & trt 0 & trt 1 & trt 0 \\
  \hline
$400$ &I& $(-0.002,0.022)$ & $(-0.003,0.044)$ & $(0.003,0.029)$ & $(0.007,0.045)$ \\
   & A & $(-0.002,0.022)$ & $(-0.003,0.043)$ & $(0.003,0.028)$ & $(0.006,0.044)$ \\
  $600$ & I & $(0.001,0.044)$ & $(-0.006,0.051)$ & $(0.013,0.055)$ & $(0.010,0.054)$ \\
   & A & $(0.003,0.042)$ & $(-0.007,0.052)$ & $(0.011,0.053)$ & $(0.008,0.054)$ \\
  $800$ & I & $(0.008,0.057)$ & $(-0.001,0.068)$ & $(0.014,0.066)$ & $(0.009,0.069)$ \\
   & A & $(0.007,0.056)$ & $(-0.003,0.067)$ & $(0.012,0.064)$ & $(0.008,0.069)$ \\
  $1000$ & I & $(0.006,0.059)$ & $(-0.005,0.080)$ & $(0.010,0.076)$ & $(0.014,0.083)$ \\
   & A & $(0.004,0.058)$ & $(-0.006,0.079)$ & $(0.010,0.072)$ & $(0.010,0.082)$ \\
   \hline
\end{tabular}

$^\dagger$I denotes the IPSWKME and A denotes the AIPSWKME; trt represents treatment; Norm CI denotes the confidence interval obtained using normal approximation based on asymptotic results; Boot CI denotes the confidence interval obtained using 500 bootstraps.
\label{tab:compare}
\end{table}

\end{document}